\documentclass[11pt,a4paper]{article}

\usepackage{amsmath,amsfonts,amssymb,latexsym,graphics,epsfig,url}
\usepackage{xcolor}
\usepackage{amsthm}
\usepackage[english]{babel}
\usepackage{mathdots}
\usepackage{graphicx}
\usepackage{soul}
\usepackage[utf8]{inputenc}
\usepackage{diagbox}
\usepackage[makeroom]{cancel}

\usepackage{float}              
\newtheorem{theorem}{Theorem}[section]
\newtheorem{proposition}[theorem]{Proposition}
\newtheorem{lemma}[theorem]{Lemma}
\newtheorem{corollary}[theorem]{Corollary}

\newtheorem{definition}[theorem]{Definition}

\DeclareMathOperator{\Aut}{Aut}

\DeclareMathOperator{\circu}{circ}

\DeclareMathOperator{\diag}{diag}
\DeclareMathOperator{\dist}{dist}
\DeclareMathOperator{\rank}{rank}
\DeclareMathOperator{\tr}{tr}
\DeclareMathOperator{\spec}{sp}

\DeclareMathOperator{\Irep}{Irep}

\def\Q{\ns Q}
\def\R{\ns R}
\def\Z{\ns Z}

\def\f{\mbox{\boldmath $f$}}

\def\u{\mbox{\boldmath $u$}}
\def\x{\mbox{\boldmath $x$}}
\def\y{\mbox{\boldmath $y$}}

\def\vecv{\mbox{\boldmath $v$}}

\def\w{\mbox{\boldmath $w$}}

\def\vec0{\mbox{\boldmath $0$}}

\def\A{\mbox{\boldmath $A$}}
\def\B{\mbox{\boldmath $B$}}
\def\C{\mbox{\boldmath $C$}}

\def\D{\mbox{\boldmath $D$}}
\def\F{\mbox{\boldmath $F$}}

\def\I{\mbox{\boldmath $I$}}
\def\J{\mbox{\boldmath $J$}}

\def\L{\mbox{\boldmath $L$}}
\def\M{\mbox{\boldmath $M$}}
\def\O{\mbox{\boldmath $O$}}

\def\X{\mbox{\boldmath $X$}}
\def\Y{\mbox{\boldmath $Y$}}
\def\Z{\ns{Z}}
\def\I{\mbox{\boldmath $I$}}
\def\J{\mbox{\boldmath $J$}}

\def\P{\mbox{\boldmath $P$}}
\def\Q{\mbox{\boldmath $Q$}}
\def\R{\mbox{\boldmath $R$}}
\def\S{\mbox{\boldmath $S$}}
\def\ZZ{\mbox{\boldmath $Z$}}
\def\1{\mbox{\boldmath $1$}}

\def\Re{\mathbb R}
\def\Z{\mathbb Z}

\newcommand\restr[2]{\ensuremath{\left.#1\right|_{#2}}}

\begin{document}
	
\title{On the spectra 
of token graphs of cycles and other graphs
\thanks{This research has been supported by
AGAUR from the Catalan Government under project 2021SGR00434 and MICINN from the Spanish Government under project PID2020-115442RB-I00.
The research of M. A. Fiol was also supported by a grant from the  Universitat Polit\`ecnica de Catalunya with references AGRUPS-2022 and AGRUPS-2023.}
}
	\author{M. A. Reyes$^a$, C. Dalf\'o$^a$, M. A. Fiol$^b$, and A. Messegu\'e$^a$\\
		\\
		{\small $^a$Dept. de Matem\`atica, Universitat de Lleida, Lleida/Igualada, Catalonia}\\
		{\small {\tt \{monicaandrea.reyes,cristina.dalfo,visitant.arnau.messegue\}@udl.cat}}\\
		{\small $^{b}$Dept. de Matem\`atiques, Universitat Polit\`ecnica de Catalunya, Barcelona, Catalonia} \\
		{\small Barcelona Graduate School of Mathematics} \\
		{\small  Institut de Matem\`atiques de la UPC-BarcelonaTech (IMTech)}\\
		{\small {\tt miguel.angel.fiol@upc.edu} }\\
	}

\date{}
\maketitle
	
\begin{abstract}
The $k$-token graph $F_k(G)$ of a graph $G$ is the graph whose vertices are the $k$-subsets of vertices from $G$, two of which being adjacent whenever their symmetric difference is a pair of adjacent vertices in $G$.
It is a known result that the algebraic connectivity  (or second Laplacian eigenvalue) of $F_k(G)$ equals the algebraic connectivity of $G$. 
 In this paper, we first give results that relate the algebraic connectivities of a token graph and the same graph after removing a vertex.
Then, 
we prove the result on the algebraic connectivity of 2-token graphs for two infinite families: the odd graphs $O_r$ for all $r$, and the multipartite complete graphs $K_{n_1,n_2,\ldots,n_r}$ for all $n_1,n_2,\ldots,n_r$
 In the case of cycles, we present a new method that allows us to compute the whole spectrum of $F_2(C_n)$. This method also allows us to obtain closed formulas that give asymptotically exact approximations for most of the eigenvalues of $F_2(\textit{}C_n)$.
\end{abstract}
	

	

\section{Introduction}
\label{sec:-1https://www.overleaf.com/project/635919a7a2694a23b29b9dab}

 Let $G=(V,E)$ be a simple graph with vertex set $V=V(G)=\{1,2,\ldots,n\}$ and edge set $E=E(G)$. By convenience, we consider every edge $e=\{u,v\}$ constituted by two opposite arcs $(u,v)$ and $(v,u)$. Let $N(u)$ denote the set of vertices adjacent to $u\in V$, so that the minimum degree of $G$ is $\delta(G)=\min_{u\in V}|N(u)|$. For a given integer $k$ such that $1\leq k \leq n$, the \textit{$k$-token graph} $F_k(G)$ of $G$ is the graph whose vertex set $V (F_k(G))$ consists of the ${n \choose k}$ $k$-subsets of vertices of $G$, and two vertices $A$ and $B$ of $F_k(G)$ are adjacent if and only if their symmetric difference $A \bigtriangleup B$ is a pair $\{a,b\}$ such that $a\in A$, $b\in B$, and $(a,b)\in E(G)$.
 Then, if $G$ has $n$ vertices and $m$ edges, $F_k(G)$ has ${n \choose k}$ vertices and ${n-2 \choose k-1}m$ edges. (Indeed, for each edge of $G$, there are 
${n-2 \choose k-1}$ edges of $F_k(G)$.)
 We also use the notation $\{a,b\}$, with $a,b\in V$, for a vertex of a 2-token graph. Moreover, we use $ab$ for the same vertex in the figures. The naming `token graph'
comes from an observation in
	Fabila-Monroy, Flores-Pe\~{n}aloza,  Huemer,  Hurtado,  Urrutia, and  Wood \cite{ffhhuw12}, that vertices of $F_k(G)$ correspond to configurations
	of $k$ indistinguishable tokens placed at distinct vertices of $G$, where
	two configurations are adjacent whenever one configuration can be reached
	from the other by moving one token along an edge from its current position
	to an unoccupied vertex. 
The $k$-token graphs are also called symmetric $k$-th power of graphs in Audenaert, Godsil, Royle, and Rudolph \cite{agrr07}, and $k$-tuple vertex graphs in Alavi, Lick, and Liu \cite{all02}.
In Figures \ref{F_2(C_9)}, \ref{F_2(C_{10})}, and  \ref{F_2(C_8)}, we show the 2-token graphs of cycles $C_9$, $C_{10}$,  and $C_8$, respectively.
Note that if $k=1$, then $F_1(G)\cong G$; and if $G$ is the complete graph $K_n$, then $F_k(K_n)\cong J(n,k)$, where $J(n,k)$ denotes the Johnson graph~\cite{ffhhuw12}, which is distance-transitive (and, hence, distance-regular).
Moreover, if $G$ is bipartite, so it is $F_k(G)$ for any $k=1,\ldots,|V|-1$.
	
Token graphs have some applications in physics. For instance, a relationship between token graphs and the exchange of Hamiltonian operators in
quantum mechanics is given in Audenaert, Godsil, Royle, and Rudolph \cite{agrr07}.
Our interest in the study of token graphs is motivated by some of their applications in 
mathematics and computer science: Analysis of complex networks, coding theory, combinatorial designs (by means of Johnson graphs),
algebraic graph theory, enumerative combinatorics, the study of symmetric functions, etc.
	
Recently, it was conjectured by Dalf\'o, Duque, Fabila-Monroy, Fiol, Huemer, Trujillo-Negrete, and Zaragoza Mart\'{\i}nez \cite{ddffhtz21} that the algebraic connectivity of $F_k(G)$ equals the algebraic connectivity of $G$. 
After submitting the first version of this paper, the authors learned (from Fabila-Monroy \cite{f23}) that this conjecture was already known as the {\em  Aldous' spectral gap conjecture},  and it was proved in 2010 by Caputo, Ligget, and Richthammer in \cite{clr10}.
 Moreover,  Ouyang  \cite{o19} and Lew \cite{l23}
 also mentioned that this conjecture was actually solved.
 Moreover, Cesi \cite{c16} provided a simpler  proof of the so-called `octopus inequality', which is one of the main ingredients to prove the Aldous’
conjecture.
 These results were obtained 
 in completely different contexts and using distinct techniques. More precisely, they used the theory of continuous Markov chains of random walks and the so-called `interchange process'. 
 In this paper, we present an algebraic approach to this problem based on voltage graphs, and we give a new method
 that can be of interest 
 giving an alternative proof.

This paper is structured as follows. In Section \ref{resultats-coneguts}, we give the preliminaries and background together with the known results. Moreover, there are the concepts of quotient graph (or digraph) and lift graph (or digraph). In Section \ref{resultats-nous}, we give results on the algebraic connectivity of a $k$-token graph when we remove one of its vertices. Besides, given 
that the algebraic connectivities of a $k$-token graph and its original graph are the same, we provide an algebraic proof of this result 
when $k=2$ for two infinite families of graphs: the odd graphs and the multipartite complete graphs.
Section \ref{sec:2-tokens} deals with the $2$-token graph of a cycle $C_n$ and gives an efficient method to compute the whole spectrum of $F_2(C_n)$ for any $n$ by using the theory of lift graphs and a new method called over-lifts.
Finally, in the last section, we give some closed formulas that provide asymptotic approximations of the eigenvalues of $F_2(C_n)$.

	
\section{Preliminaries and background}
\label{resultats-coneguts}

\subsection{Some notation and basic facts}

The notation of this paper is as follows: 
$\A=\A(G)$ is the adjacency matrix of the graph $G$, $\L=\L(G)$ the Laplacian matrix of the graph $G$, $\P$ a permutation matrix, $\pi=V_1\cup V_2\cup \cdots \cup V_r$ a (regular or not) partition of the vertex set, $G/\pi$ a quotient graph over $\pi$, $\A(G/\pi)$ and $\L(G/\pi)$ the adjacency and Laplacian matrices of $G/\pi$, and $\S$ the characteristic matrix of the partition $\pi$.

The transpose of a matrix $\M$ is denoted by $\M^\top$, the
identity matrix by $\I$, the all-$1$ vector $(1,\ldots, 1)^{\top}$ by $\1$, the all-$1$ (universal) matrix  by $\J$, and the all-$0$ vector and all-$0$ matrix by $\vec0$
and $\O$, respectively.
Let $[n]:=\{1,\ldots,n\}$ and ${[n]\choose k}$ denote the set of $k$-subsets of $[n]$, which is the set of vertices of the $k$-token graph.
	
For our purpose, it is convenient to denote by $W_n$ the set of all column vectors  $\vecv$ with $n$ entries such that $\vecv^{\top }\1 = 0$.
Recall that any square matrix $\M$ with all zero row sums has an eigenvalue $0$ with corresponding eigenvector $\1$. Then, given a graph $G=(V,E)$ of order $n$, we say that a vector $\vecv\in \mathbb{R}^n$ is an \textit{embedding} of $G$ if $\vecv\in W_n$. Note that if $\vecv$ is a $\lambda$-eigenvector of $G$, with $\lambda>0$, then it is an embedding of $G$.
	
When $\M=\L(G)$, the Laplacian matrix of a graph $G$, the matrix is positive semidefinite, with eigenvalues $(0=)\lambda_1\le \lambda_2\le \cdots \le \lambda_n$. Its second smallest eigenvalue $\lambda_2$ is known as the {\em algebraic connectivity} of  $G$ (see Fiedler \cite{fi73}), and we denote it by $\alpha(G)$.
The spectral radius $\lambda_{\max}(G)=\lambda_n$ satisfies several lower and upper bounds (see Patra and Sahoo \cite{ps17} for a survey).

For a graph $G$ with Laplacian matrix $\L(G)$ and an embedding $\vecv$ of $G$, let
	$$
	\lambda_G(\vecv):=\frac{\vecv^{\top}\L(G)\vecv}{{\vecv}^{\top}\vecv}=\frac{\sum\limits_{(i,j)\in E}[\vecv(i)-\vecv(j)]^2}{\sum\limits_{i\in V}\vecv^2(i)},
	$$
	where
	$\vecv(i)$ denotes the entry of $\vecv$ corresponding to the vertex $i\in V(G)$.
	The $\lambda_G(\vecv)$ value is known as
	the {\em Rayleigh quotient}.
	If $\vecv$ is an eigenvector of $G$, then its corresponding eigenvalue is $\lambda(\vecv)$.
	Moreover, for an embedding $\vecv$ of $G$, we have
	\begin{equation}
		\label{bound-lambda(v)}
		\alpha(G)\le \lambda_G(\vecv),
	\end{equation}
	and we have equality  when $\vecv$ is an $\alpha(G)$-eigenvector of $G$.

In this paper, we first give some results about the algebraic connectivity of a token graph. Besides, we provide results about the spectrum of a token graph of a cycle graph when we deal with 2 tokens. 
This study was initiated by Dalf\'o,  Duque, Fabila-Monroy,  Fiol, Huemer,  Trujillo-Negrete, and  Zaragoza Mart\'{\i}nez in \cite{ddffhtz21}. One of their results is the following. 
 
 Given some integers $n$ and $k$ (with $k\in [n]$), we define the $(n;k)$-\emph{binomial matrix} $\B$. This is a ${n \choose k}\times n$ matrix whose rows are the characteristic vectors of the $k$-subsets of $[n]$ in a given order. Thus, if the $i$-th $k$-subset is $A$, then
$$
(\B)_{ij}=
\left\lbrace
\begin{array}{ll}
	1 & \mbox{if } j\in A,\\
	0 & \mbox{otherwise.}
\end{array}
\right.
$$

 \begin{lemma}[\cite{ddffhtz21}]
		\label{coro:LkL1}
		Let $G$ be a graph on $n$ vertices. For some integers $h,k$ such that $1\le h<k\le \frac{n}{2}$,  let $F_h=F_h(G)$ and $F_k=F_k(G)$ be its $h$- and $k$-token graphs with respective Laplacian matrices $\L_h$ and $\L_k$.
		Then, the following holds:
		\begin{itemize}
			\item[$(i)$]
			If $\vecv$ is a $\lambda$-eigenvector of $\L_h$, then $\B\vecv$ is a $\lambda$-eigenvector of $\L_k$.
			Thus, the Laplacian spectrum (eigenvalues and their multiplicities) of $\L_h$ is contained in the Laplacian spectrum of $\L_k$.
			\item[$(ii)$]
			If $\u$ is a $\lambda$-eigenvector of $\L_k$ such that $\B^{\top}\u\neq \vec0$, then $\B^{\top}\u$
			is a $\lambda$-eigenvector of $\L_h$.
		\end{itemize}
\end{lemma}
From the inclusion property of the successive spectra in $(i)$, we have:
\begin{equation}
\label{eq:non-increasing}
\alpha(G)\ge \alpha(F_2(G))\ge \alpha(F_3(G))\ge\cdots\ge \alpha(F_{\lfloor n/2\rfloor}(G)).
\end{equation}
Recall that Caputo, Liggett, and Richthammer \cite{clr10} proved that all these inequalities actually are equalities.

	
	
	In our context of token graphs, it was shown in \cite{ddffhtz21} and by Dalf\'o and Fiol in \cite{df22} that the conjecture (now, a result) holds for the following infinite families of graphs.
\begin{theorem}[\cite{ddffhtz21,df22}]
		\label{theo:alg-connec-antic}
		For each of the following classes of graphs, the algebraic connectivity of a token graph $F_k(G)$ equals the algebraic connectivity of $G$.
		\begin{itemize}
			\item[$(i)$]
			Let $G=K_n$ be the complete graph on $n$ vertices. Then,
			$\alpha(F_k(G))=\alpha(G)=n$ for every $n$ and $k=1,\ldots,n-1$.
			\item[$(ii)$]
			Let $G= K_{n_1,n_2}$ be the complete bipartite graph on $n=n_1+n_2$ vertices, with $n_1\le n_2$. Then, $\alpha(F_k(G))=\alpha(G)=n_1$ for every $n_1,n_2$ and $k=1,\ldots,n-1$.
			\item[$(iii)$]
			Let $T_n$ be a tree on $n$ vertices. Then,
			$\alpha(F_k(T_n))=\alpha(T_n)$ for every $n$ and $k=1,\ldots,n-1$.
			\item[$(iv)$]
			Let $G$ be a graph  such that $\alpha(F_k(G))=\alpha(G)$. Let $T_G$ be a graph obtained from $G$ by attaching a (possibly empty) $u$-rooted tree $T(u)$ to each vertex $u$ of $G$.
   Then,
			$\alpha(F_k(T_G))=\alpha(T_G)$.
		\end{itemize}
	\end{theorem}
 
All these results were obtained by induction on $n$, and using that for some vertex $i$, $\alpha(G\setminus i)\ge \alpha(G)$. See Table \ref{tab:Data1}, where we show two particular cases of $(iii)$. Namely, the star graph $S_n$ and the path graph $P_n$, both on $n$ vertices.
These two cases, together with the complete bipartite graph $K_{n_1,n_2}$ were proved in \cite[Th. 7.2]{ddffhtz21}.  
The fact that, for any tree $T_n$, we have $\alpha(T_n\setminus i)\geq \alpha(T_n)$ can be proved by using interlacing (see Bunch, Nielsen, and Sorensen \cite{bns78}). An alternative proof using Fiedler vectors (eigenvectors with their eigenvalue equal to the algebraic connectivity) was given in Dalf\'o and Fiol \cite{df22}.

\begin{table}[H]
    \centering
    \begin{tabular}{|c|c|c|c|c|}
    \hline
     Graph $G$& \bf{ $\alpha(G)$}& Vertex $i$  & \textbf{$\alpha(G\setminus i)$}  \\
     \hline\hline
    $P_n$ & $2(1-\cos(\frac{\pi}{n}))$& a leaf  & $2(1-\cos(\frac{\pi}{n-1}))$  \\
     $S_n$ & $1$ & a leaf & $1$  \\
     $K_{n_1,n_2}$ ($n_1 < n_2$) & $n_1$ &  $i \in V_2$, $n_2 = \vert V_2\vert$   & $n_1$  \\
     $T_n$& (see \cite{df22}) & a leaf  & (see \cite{df22})  \\
     \hline
    \end{tabular}
     \caption{Some graphs with a vertex $i$ such that $\alpha(G) \leq \alpha (G\setminus i)$.}
     \label{tab:Data1}
\end{table}






\subsection{Regular partitions and their spectra}
\label{sec:reg-part}

Let $G=(V,E)$ be a graph with vertex set $V=V(G)$ and Laplacian matrix $\L$. A partition $\pi=(V_1,\ldots,
V_r)$ of $V$ is called {\em regular} (or {\em equitable})
whenever, for any $i,j=1,\ldots,r$, the {\em intersection numbers}
$b_{ij}(u)=|N(u)\cap V_j|$, where $u\in V_i$, do not depend on the vertex $u$ but only on the subsets (usually called {\em classes} or {\em
cells}) $V_i$ and $V_j$. In this case, such numbers are simply written as $b_{ij}$, and the $r\times r$ matrix $\Q_L=\L(G/\pi)$ with entries
$$
(\Q_L)_{ij}=\left\{
\begin{array}{cl}
\!\!\!\!\!\!-b_{ij} & \mbox{if } i\neq j,\\[.2cm]
b_{ii}-\displaystyle \sum_{h=1}^r b_{ih} & \mbox{if } i=j , 
\end{array}
\right.
$$
is referred to as the \textit{quotient Laplacian matrix} of $G$ with respect to $\pi$. This is also represented by the \textit{quotient (weighted) directed graph} $G/\pi$ (associated with the partition $\pi$),
with vertices representing the $r$ cells, and there is an arc with weight $b_{ij}$ from vertex $V_i$ to vertex $V_j$ if and only if $b_{ij}\neq 0$.
Of course, if $b_{ii}>0$, for some $i=1,\ldots,r$, the quotient graph (or digraph) $G/\pi$ has loops.
Given a partition $\pi$ of $V$ with $r$ cells (or partition sets), let $\S$ be the {\em characteristic matrix} of $\pi$, that is, the $n\times r$ matrix whose columns are the characteristic vectors of the cells of $\pi$. Then, as in the case of the adjacency matrix (see, for instance, Godsil and Royle \cite{gr01}),  $\pi$ is a regular partition if and only if $\L\S=\S\Q_L$. In the case of the Laplacian matrix, it follows that 
$$
\Q_L=(\S^{\top}\S)^{-1}\S^{\top}\L\S,
$$
and the characteristic polynomial of $\Q_L$ divides the characteristic polynomial of $\L$.
Thus, a part of the spectrum of $\L$ can be determined by the spectrum of the (usually much smaller) matrix $\Q_L$.
Moreover, if the graph $G$ is bipartite, the maximum eigenvalue of its Laplacian quotient matrix $\Q_L$ equals the spectral radius of $\L$. The reason is that, in bipartite graphs, the Laplacian matrix's characteristic polynomial is equal to the signless Laplacian $\L^+$ (see, for instance, Grone, Merris, and  Sunder \cite{gms90}). Then, the same holds for the quotient Laplacian $\Q_L$ and quotient signless Laplacian $\Q_L^+$  matrices. Thus, since each eigenvector of $\Q_L^+$ gives rise to an eigenvector of  $\L^+$, the spectral radius of $\L$ corresponds to the eigenvalue of the Perron vector of $\Q_L^+$ or maximum eigenvalue of $\Q_L$.
For more information about quotient (Laplacian) matrices, see 
Dalf\'o, Fiol, Pavl\'ikov\'a, and \v{S}ir\'an \cite{df22b}.

\subsection{Lift graphs and their spectra}
\label{sec:sp}
Let $\cal{G}$ be a group. An ({\em ordinary\/}) {\em voltage assignment} on the (di)graph (that is, graph or digraph) $G=(V,E)$  is a mapping $\beta: E\to \cal{G}$ with the property that $\beta(a^-)=(\beta(a^+))^{-1}$ for every arc $a\in E$. Thus, a voltage assigns an element $g\in \cal{G}$ to each arc of the (di)graph so that a pair of mutually reverse arcs $a^+$ and $a^{-}$, forming an undirected edge, receive mutually inverse elements $g$ and $g^{-1}$. The (di)graph $G$ and the voltage assignment $\beta$ determine a new (di)graph $G^{\beta}$, called the {\em lift} of $G$, which is defined as follows. The vertex and arc sets of the lift are simply the Cartesian products $V^{\beta}=V\times \cal{G}$ and $E^{\beta}=E\times \cal{G}$, respectively. Moreover, for every arc $a\in E$ from a vertex $u$ to a vertex $v$ for $u,v\in V$ (possibly, $u=v$) in $G$, and for every element $g\in \cal{G}$, there is an arc $(a,g)\in E^{\beta}$ from the vertex $(u,g)\in V^{\beta}$ to the vertex $(v,g\beta(a))\in V^{\beta}$.

Let $G=(V,E)$ be a connected graph on $n$ vertices (with loops and multiple edges allowed) and with Laplacian matrix $\L$. Let $\beta$
be a voltage assignment on the arc set $E$  in a  group
$\cal{G}$ with identity element $e$.
Now we show that the spectrum of the Laplacian matrix of the lift $G^{\beta}$  may be computed.
To this end, the key idea is to define the so-called {\em Laplacian base matrix} properly as follows.

\begin{definition}
\label{B(U)}
To the pair $(G,\beta)$, we assign the $n\times n$ {\em Laplacian base matrix} $\B(\L)$ defined by 
$$
\B(\L)=-\B(\A)+\B(\D),
$$
where the matrices $\B(\A)$ and $\B(\D)$ have entries as follows:
\begin{itemize}
\item
$\B(\A)_{uv}=\beta(a_1)+\cdots +\beta(a_j)$ if $a_1,\ldots,a_j$ is the set of all the arcs of $G$ from $u$ to $v$, not excluding the case $u=v$, and $\B(\A)_{uv}=0$ if $(u,v)\not\in E$.
\item
$\B(\D)_{uu}=\deg(u)\cdot e$, and $\B(\D)_{uv}=0$ if $u\neq v$.
\end{itemize}
\end{definition}

Let $\rho \in \Irep(\cal{G})$ be a unitary irreducible representation of $\cal{G}$ of dimension $d_\rho=\dim(\rho)$. Given a graph $G$ on $n$ vertices, the assignment $\beta$ in $G$, and the Laplacian base matrix $\B=\B(\L)$, let $\rho(\B)$ be the $d_\rho n\times d_\rho n$ matrix obtained from $\B$ by replacing every nonzero entry $(\B)_{u,v} \in \mathbb{C}[\cal{G}]$ as above by the $d_\rho\times d_\rho$ matrix $\rho(\B_{u,v})$. That is, each element $g$ of the group is replaced by $\rho(g)$, and
the zero entries of $\B$ are changed to all-zero $d_\rho\times d_\rho$ matrices. We refer to $\rho(\B)$ as the {\em $\rho$-image} of the Laplacian base matrix $\B$.
For every $\rho\in \Irep(\cal{G})$, we consider the $\rho$-image $\rho(\B)$ of the Laplacian base matrix $\B$, and we let $\spec(\rho(\B))$ denote the spectrum of $\rho(\B)$, that is, the multiset of all the $d_\rho n$ eigenvalues of the matrix $\rho(\B)$. Finally, the notation $d_\rho\cdot\spec(\rho(\B))$ denotes the multiset 
obtained by taking each of the $d_\rho n$ entries of the spectrum $\spec(\rho(\B))$ exactly $d_\rho$ times. In particular, if $(d_\rho)=0$, we take $d_\rho\cdot\spec(\rho(B))=\emptyset$.

With all these notations, we can state the following result from 
Dalf\'o, Fiol, Pavl\'ikov\'a, and \v{S}ir\'a\v{n} \cite{df22b}, in which, to our knowledge, the theory of lift graphs was applied to the study of the spectrum of token graphs for the first time. This result  
allows us to compute the spectrum of a (regular) lifted (di)graph from its associated matrix and the irreducible representations
of its corresponding group. 

\begin{theorem}[\cite{dfs19}]
\label{theo-sp}
Let $G=(V,E)$ be a base (di)graph on $n$ vertices, with a voltage assignment $\beta$ in a group $\cal{G}$.
For every irreducible representation $\rho\in \Irep (\cal{G})$, let $\rho(\B)$ be the complex matrix whose entries are given by $\rho(\B_{u,v})$.
Then,
$$
\spec G^{\beta} = \bigcup_{\rho\in \Irep(\cal G)}d_\rho\cdot\spec(\rho(\B)).
$$
\end{theorem}

In Section \ref{sec:2-tokens}, we use this result in the case when the group $\cal{G}$ is cyclic. Then, if $g$ is a generator of $\cal{G}$, with order $n$, the faithful representation $\rho$ such that $\rho(g^r)=\zeta^r$, with $\zeta=e^{i\frac{2\pi}{n}}$, has dimension $1$. Then, we consider the Laplacian base matrix with each entry being a polynomial in $z$ with integer coefficients 
and represent such a `polynomial matrix' by $\B(z)$. Thus, Theorem \ref{theo-sp}
gives
$$
\spec G^{\beta} = \bigcup_{z\in R(n)}\spec(\B(z)),
$$
where  $R(n)$ is the set of all $n$-th roots of unity. 
A simple property of the polynomial matrix is that $\B(1)$ is the quotient matrix of a regular partition of the lift graph $G^{\beta}$. For more information on lift graphs, see Dalf\'o, Fiol, Miller, Ryan, and \v{S}ir\'a\v{n} \cite{dfmrs17}.


\section{The algebraic connectivity of $F_k(G)$}
\label{resultats-nous}

In this section, we give some results on the algebraic connectivity of token graphs. We begin with a known lemma and continue with a new result, from which some of the algebraic connectivity results are derived.

Let $G$ be a graph with $k$-token graph $F_k(G)$.
For a vertex $a\in V(G)$, let $S_a:=\{A\in V(F_k(G)):a\in A\}$ and $S'_a:=\{B\in V( F_k(G)): a\not\in B\}$.
Let $H_a$ and $H'_a$ be the subgraphs of $F_k(G)$ induced by $S_a$ and $S'_a$, respectively.
Note that $H_a\cong F_{k-1}(G\setminus \{a\})$ and $H'_a\cong F_k(G\setminus \{a\})$.
\begin{lemma}[\cite{ddffhtz21}]
\label{lem:embedding}
Given a vertex $a\in G$ and an eigenvector $\vecv$ of $F_k(G)$ such that $\B^{\top}\vecv=\vec0$, let
$$
\w_a:=\restr{\vecv}{S_a} \mbox{ and } \quad \w'_a:=\restr{\vecv}{S'_a}.
$$
Then, $\w_a$ and $\w'_a$ are embeddings of $H_a$ and $H'_a$, respectively.
\end{lemma}

\begin{lemma}
\label{lem:basic}
 Let $G=(V,E)$ be a graph with token graph $F_k=F_k(G)$ for some integer $k\ge 2$. 
 Let $\xi(G^-)=\min_{i\in V}\alpha(F_{k-1}(G\setminus i))$. If $\alpha(F_k(G))<\alpha(G)$, then the following statements hold.
 \begin{itemize}
 \item[$(i)$] 
 $\alpha(F_k(G))\ge \frac{k}{k-1}\xi(G^-)$,
 \item[$(ii)$] 
 $\alpha(F_k(G))\ge \alpha(F_k(G\setminus i))$ for $i\in V$.
 \end{itemize}
\end{lemma}
\begin{proof}
If $\alpha(F_k(G))<\alpha(G)$, we know that the eigenvector $\vecv$ of 
$\alpha(F_k(G))$ must satisfy $\B^{\top}\vecv=\vec0$. 
Let $\|\vecv\|=1$.
Given a vertex $i\in V$, let $S_i:=\{A\in V(F_k):i\in A \}$ and   
$S'_i:=\{B\in V(F_k):i \not\in B \}$.
Let $H_i\cong F_{k-1}(G\setminus i)$ and  $H'_i\cong F_{k}(G\setminus i)$ be the subgraphs of $F_k(G)$ induced by $S_i$ and $S'_i$, respectively.
Let $\w_i=\vecv|_{S_i}$ and $\w'_i=\vecv|_{S'_i}$. 
 From $\B^{\top}\vecv=\vec0$, from Lemma \ref{lem:embedding}, 
 $\w_i$ and $\w'_i$ are embeddings of $H_i$ and $H'_i$, respectively. Then, by \eqref{bound-lambda(v)}, we have that their Rayleigh quotients satisfy
\begin{align}
    \lambda(\w_i)& =\frac{\sum\limits_{(A,B)\in E(H_i)} [\w_i(A)-\w_i(B)]^2}{\sum\limits_{A\in V(H_i)}\w_i(A)^2}\ge \alpha(F_{k-1}(G\setminus i))\ge \alpha(F_{k}(G\setminus i)),
    \label{eq:Rlambda(w)}\\
     \lambda(\w'_i)& =\frac{\sum\limits_{(A,B)\in E(H'_i)} [\w'_i(A)-\w'_i(B)]^2}{\sum\limits_{A\in V(H'_i)}\w'_i(A)^2}\ge \alpha(F_{k}(G\setminus i)),
     \label{eq:Rlambda(w')}
\end{align}
where, in the last inequality of \eqref{eq:Rlambda(w)}, we applied \eqref{eq:non-increasing} with $G\setminus i$.
In order to prove $(i)$, we use \eqref{eq:Rlambda(w)}:
\begin{align*}
\alpha(F_k) =\lambda(\vecv)&=\sum\limits_{(A,B)\in E(F_k)}[\vecv(A)-\vecv(B)]^2
=\frac{1}{k-1}\sum_{i=1}^n \sum\limits_{(A,B)\in E(H_i)}[\w_i(A)-\w_i(B)]^2  \\
& \ge \frac{1}{k-1}\xi(G^-)\sum_{i=1}^n\sum\limits_{A\in V(H_i)}\vecv(A)^2=
\frac{k}{k-1}\xi(G^-),
\end{align*}
since, in the double summatory, each edge $(A,B)$ of $F_k$ is considered $k-1$ times, whereas for the last equality, each vertex $A$ of $F_k$ is considered $k$ times.
Now, we prove $(ii)$ by using \eqref{eq:Rlambda(w)} and \eqref{eq:Rlambda(w')}. Since 
$V(H_i)\cup V(H'_i)=V(F_k)$, and $\|\vecv\|=1$, we have:
\begin{align*}
\alpha(F_k) =\lambda(\vecv)&=\sum\limits_{(A,B)\in E(F_k)}[\vecv(A)-\vecv(B)]^2\\
 &\ge \sum\limits_{(A,B)\in E(H_i)}[\w_i(A)-\w_i(B)]^2
 +\sum\limits_{(A,B)\in E(H'_i)}[\w'_i(A)-\w'_i(B)]^2\\
 & \ge \alpha(F_k(G\setminus i))\left[\sum\limits_{A\in V(H_i)}\w_i(A)^2+\sum\limits_{B\in V(H'_i)}\w'_i(B)^2\right]\\
 & \ge \alpha(F_k(G\setminus i))\left[\sum\limits_{A\in V(H_i)}\vecv(A)^2+\sum\limits_{B\in V(H'_i)}\vecv(B)^2\right]=\alpha(F_k(G\setminus i)).
 \label{lem:basic}
\end{align*}
This completes the proof.
\end{proof}
\begin{theorem}
Let $G$ be a graph with a vertex $i$ such that $\alpha(G)\ge \alpha(G\setminus i)$. Then, $\alpha(F_k(G))\ge \alpha(F_k(G\setminus i)).$ 
\end{theorem}
\begin{proof}
If $\alpha(F_k(G))<\alpha(G)$,  Lemma \ref{lem:basic}$(ii)$ gives $\alpha(F_k(G))\ge \alpha(F_k(G\setminus i))$. Then, if $\alpha(F_k(G))< \alpha(F_k(G\setminus i))$, we must have,
using \eqref{eq:non-increasing},
$$
\alpha(F_k(G))\ge\alpha(G)\ \Rightarrow\ \alpha(F_k(G))=\alpha(G)\ \Rightarrow\ \alpha(G)<\alpha(F_k(G\setminus i))\le\alpha(G\setminus i).
$$
Thus, the result corresponds to the contrapositive statement.
\end{proof}

Now we consider the case of 2 tokens and vertex-transitive graphs.  In the following result, given a vertex $i$,  we use the parameter $\kappa(i)=\frac{\alpha(G\setminus i)}{\alpha(G)}$ introduced by Kirkland in \cite{k10}. When $G$ is vertex-transitive, we denote $\kappa(G)=\kappa(i)$ for all $i$, so that $\xi(G^-)=\alpha(G\setminus i)=\kappa(G)\alpha(G)$.

\begin{theorem}
\label{th:G-vt}
Let $G=(V,E)$ be a graph on $n>3$ vertices.
\begin{itemize}
\item[$(i)$] 
If $\displaystyle\min_{i\in V}\kappa(i)\ge \textstyle\frac{1}{2}$,
then $\alpha(F_2(G))=\alpha(G)$.
\item[$(ii)$] 
If $G$ is vertex-transitive and $\kappa(G)\ge \frac{1}{2}$,
then $\alpha(F_2(G))=\alpha(G)$.
\end{itemize}
\end{theorem}

\begin{proof}
Let $\vecv$ be an eigenvector of $F_2=F_2(G)$ with eigenvalue $\alpha(F_2)$ and norm $\|\vecv\|=1$. If $\B^{\top}\vecv\neq \vec0$, by Lemma \ref{coro:LkL1}$(ii)$, $\alpha(F_2)=\alpha(G)$, and the results hold. Thus, we can assume that $\B^{\top}\vecv=\vec0$
(this could only occur if  $F_2\neq G$).
To prove $(i)$, we use Lemma \ref{lem:basic}$(i)$ with $k=2$ that yields
$$
\alpha(F_2)\ge 2\cdot \min_{i\in V}\alpha(G\setminus i)=2\alpha(G)\cdot\min_{i\in V}\kappa(i).
$$
Thus, if $\min_{i\in V}\kappa(i)\ge \frac{1}{2}$, we have $\alpha(F_2)\ge \alpha(G)$ and the claimed equality is obtained.

Now, $(ii)$ is a consequence of $(i)$ since, when $G$ is vertex-transitive, $\alpha(G\setminus i)=\alpha (G\setminus j)$ for every $i,j\in V$ and, then,
$\min_{i\in V}\kappa(i)=\kappa(G)$.
\end{proof}

In Table  \ref{tab:Data2}, there are some examples of graphs that satisfy Theorem \ref{th:G-vt}$(ii)$.

\begin{table}[H]
    \centering
    \begin{tabular}{|c|c|c|c|}
    \hline
     Graph $G$& \bf{ $\alpha(G)$} & \textbf{$\alpha(G\setminus i)$}\\
     \hline\hline
     $K_n$ & $n$&  $n-1$  \\
    Petersen & $2$& $\approx 1.26$  \\
    Heawood & $3-\sqrt{2} \approx 1.58$& $1$ \\
    Tetrahedron & $4$& $3$ \\
    Octahedron & $4$&  $3$ \\
    Hexahedron & $2$& $2(1-\cos(\frac{2\pi}{5}))\approx1.38$ \\
    Dodecahedron & $3-\sqrt{5}\approx0.76$& $\approx 0.59$ \\
    Icosahedron & $5-\sqrt{5}\approx 2.76$& $\frac{422009\pi}{605811} \approx 2.18$ \\
    Truncated tetrahedron  & $1$ & $\approx 0.53$ \\
    Prism Graph $GP(n,1)$  $(n=3,4)$ & $2$ & $\approx 1.38$  \\
    Prism Graph $GP(n,1)$ $(n>4)$ & $2\left(1-\cos\left(\frac{2\pi}{n}\right)\right)$ & $-$  \\
    Hypercube $Q_n$ & 2 & $\ge 1$\\
     \hline
    \end{tabular}
    \caption{Some vertex-transitive graphs, with $\alpha(G) \leq 2\alpha (G\setminus i)$, where $i$ is a vertex.}
    \label{tab:Data2}
\end{table}

Note that, in most cases of Table \ref{tab:Data2}, $\alpha(G)\geq2$. In fact, this condition is sufficient because Fiedler \cite{fi73} proved that, for any graph, $\alpha(G\setminus i)\ge \alpha(G)-1$ or, equivalently, when $G$ is vertex-transitive,  $\kappa(G)\ge 1-\frac{1}{\alpha(G)}$. From this, if $\alpha(G)\ge 2$, the condition $\kappa(G)\ge \frac{1}{2}$ holds, as required to have $\alpha(F_2(G))=\alpha(G)$. An alternative, much more involved, proof of this case is given in Dalf\'o, Fiol, and Messegu\'e \cite{dfm22}.

Besides the hypercube $Q_n$ cited in Table \ref{tab:Data2}, the following graphs constitute another infinite family of vertex-transitive graphs with algebraic connectivity $2$.
The odd graph $O_r$ has vertices labeled with the $(r-1)$-subsets of a $(2r-1)$-set, and two vertices are adjacent if their corresponding subsets have void intersection. Thus, $O_r$ has ${2r-1\choose r-1}$ vertices, and it is regular of degree $r$. For instance, $O_2=K_3$, and $O_3$ is the Petersen graph. 
\begin{corollary}
      Let $O_r$ be the odd graph of degree $r$. 
      Then,
$\alpha(F_2(O_r))=\alpha(O_r)$.
\end{corollary}
\begin{proof}
As the second largest (adjacency matrix) eigenvalue of $O_r$ is $r-2$ for $r\ge 3$, its algebraic connectivity is $\alpha(O_r)=r-(r-2)=2$, as claimed. Then, according to the above consequences of Theorem \ref{th:G-vt}, the algebraic connectivities of $F_2(O_r)$ and $O_r$ coincide.
\end{proof}

Another infinite family satisfying Theorem \ref{th:G-vt}$(i)$ is
the following.
\begin{corollary}
Let $G=K_{n_1,n_2,\ldots,n_r}\neq K_r$ be the multipartite complete graph on $n=\sum_{i=1}^r n_i$ vertices with $n_1\le n_2\le \cdots \le n_r$, and $r\geq3$. Then,
$\alpha(F_2(G))=\alpha(G)$.
\end{corollary}
\begin{proof}
Notice that $\overline{G}=K_{n_1}\cup \cdots \cup K_{n_r}$. Then, the maximum eigenvalue of $\overline{G}$ coincides with the maximum eigenvalue of $K_{n_r}$, which is $\lambda_{\max}(K_{n_r})=n_r$. Thus, 
$\alpha(G)=n-\lambda_{\max}(\overline{G})=n-n_r$ (see Fiedler \cite{fi73}). Thus, by taking any vertex $i\not\in V_r$  (with $|V_r|=n_r$), we have $\xi(G^-)=\alpha(G\setminus i)=n-n_r-1$. Thus, $\min_{i\in V}\kappa(i)=\frac{n-n_r-1}{n-n_r}\ge \frac{1}{2}$, since $r\geq3$ and $n\ge 4$
(recall that $G\neq K_r$). Then, Theorem \ref{th:G-vt}$(i)$ gives the result.
\end{proof}

\section{The spectrum of $F_2(C_n)$}
\label{sec:2-tokens}

In this section, we provide some results about the Laplacian spectrum (and, in particular, the algebraic connectivity) of the $k$-token $F_k(C_n)$ of the cycle $C_n$ on $n$ vertices $0,1,\ldots,n-1$, where $i$ is adjacent to $i+1\ (\textrm{mod}\, n)$. By Lemma \ref{coro:LkL1}$(i)$, we already know that $F_k(C_n)$ contains all the eigenvalues of $C_n$. Namely,
 \begin{equation}
\label{tj}
\theta_{j} =2\left(1-\cos\left(\frac{j2\pi}{n}\right)\right)=4\sin^2\left(\frac{j\pi}{n}\right),\qquad j=0,1,\ldots, n-1.
\end{equation}
Next, we begin with a general lower bound for $\alpha(F_k(C_n))$, which is given in terms of the algebraic connectivity of the path $P_{n-1}$ on $n-1$ vertices.

\begin{theorem}
\label{th:cn-vs-Pn}
Let $C_n$ be the cycle graph on $n>3$ vertices. Then,
\begin{equation}
\alpha(F_k(C_n))\ge \frac{k}{k-1}\alpha(P_{n-1}) =\frac{2k}{k-1}\left(1-\cos\left(\frac{\pi}{n-1}\right)\right)
\label{eq:bound-token-cycle}
\end{equation}
for every $n$ and $k=2,\ldots, \lfloor n/2\rfloor$.
\end{theorem}

\begin{proof}
As before, let $\vecv$ be an eigenvector of $F_k=F_k(C_n)$ with eigenvalue $\alpha(F_k)$ and norm $\|\vecv\|=1$. If $\B^{\top}\vecv\neq \vec0$, by Lemma \ref{coro:LkL1}$(ii)$, $\alpha(F_k(C_n))=\alpha(C_n)$, and the result holds since $\alpha(C_n)>2\alpha(P_{n-1})$. 
Thus, since $n\neq 3$, we can assume that $\B^{\top}\vecv=\vec0$
and apply again Lemma \ref{lem:basic}$(i)$ to get
\begin{align}
\alpha(C_n) 
   & \ge \frac{k}{k-1}\xi(C_n^-)=\frac{k}{k-1}\alpha(F_{k-1}(C_n\setminus i))=\frac{k}{k-1}\alpha(F_k(P_{n-1}))\\
   &=\frac{k}{k-1}\alpha(P_{n-1}),
\end{align}
where the last equality follows from Theorem \ref{theo:alg-connec-antic}$(iii)$.
\end{proof}
In particular, for $k=2$, \eqref{eq:bound-token-cycle} gives
that $\alpha(F_2(C_n))\ge 2\alpha(P_{n-1})$, and equality holds for $n=4$
since $\alpha(F_2(C_4))= 2\alpha(P_{3})=2$.

\begin{figure}[t]
\begin{center}
\includegraphics[width=14cm]{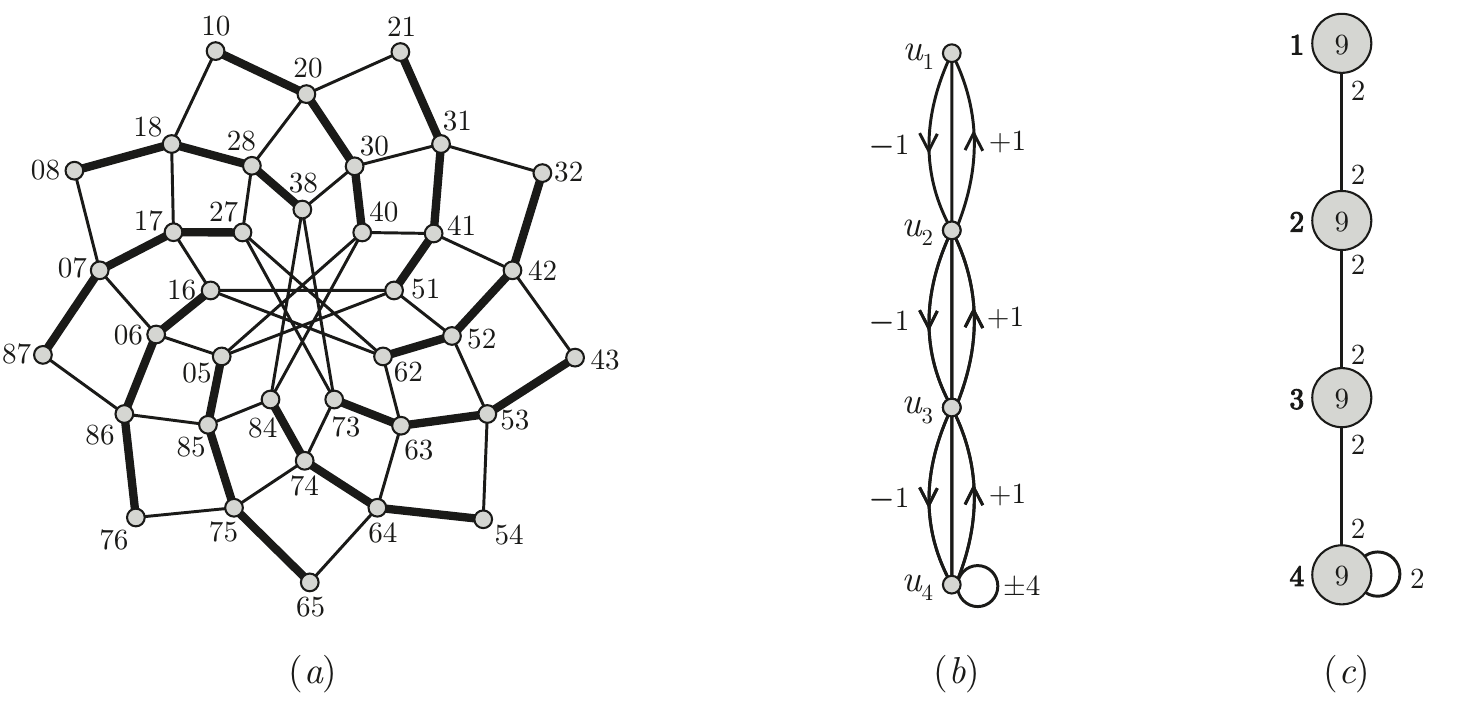}
\caption{$(a)$ The 2-token graph $F_2(C_9)$ of the cycle graph $C_9$. The thick edges correspond to each of the copies of the base graph or the quotient graph. $(b)$
Its base graph with voltages on $\mathbb{Z}_9$. $(c)$ The quotient graph of its path-shaped regular partition. In boldface, there is the numbering of the vertex classes. In class $c\in\{1,2,3,4\}$, there are the vertices $ij$ that satisfy $i-j=c\,(\textrm{mod}\, 9)$. }
\label{F_2(C_9)}
\end{center}
\end{figure}

\begin{proposition}
\label{propo:path-shaped}
\begin{itemize}
\item[$(i)$]
If $n=2\nu$ is even,
then the 2-token graph $F_2(C_n)$ has the eigenvalues
\begin{equation}
\lambda_{r} =8\sin^2 \left(\frac{r\pi}{n-1}\right),\qquad \mbox{$r=0,1,\ldots, \nu-1,$} \label{lr-even}
\end{equation}
with $\lambda_0=0$, and $\lambda_{\nu-1}=8\sin^2\left(\frac{n-2}{n-1}\frac{\pi}{2}\right)$ being the spectral radius of $F_2(C_n)$.
\item[$(ii)$]
If $n=2\nu+1$ is odd,
then the 2-token graph $F_2(C_n)$ has the eigenvalues
\begin{equation}
\lambda_{r} =8\cos^2 \left(\frac{r\pi}{n-1}\right),\qquad \mbox{$r=1,2,\ldots, \nu,$} 
\label{lr-odd}
\end{equation}
with $\lambda_{\nu}=0$, and $\lambda_{1}=8\cos^2\left(\frac{\pi}{n-1}\right)$ being a lower bound for the spectral radius of $F_2(C_n)$.
\end{itemize}
\end{proposition}
\begin{proof}
Let us see that $F_2=F_2(C_n)$ has a regular `path-shaped' partition $\pi$ with  $r=\lfloor n/2 \rfloor$ classes $V_1,V_2,\ldots,V_{r}$, where $V_i$ consists of the vertices $\{u,v\}$ such that $\dist(u,v)=i$ in $C_n$. 
\begin{itemize}
\item 
Each vertex $\{u,v\}$ in $V_1$ is adjacent to $2$ vertices
$\{u+1,v\}$ and $\{u,v+1\}$ in $V_2$ (all arithmetic is modulo $n$).
\item 
Each vertex $\{u,v\}$ in $V_i$, for $i=2,\ldots,r-1$, is adjacent to $2$ vertices $\{u-1,v\}$ and $\{u,v-1\}$ in $V_{i-1}$, and  $2$ vertices $\{u+1,v\}$ and $\{u,v+1\}$ in $V_{i+1}$.
\item
Every vertex $\{u,v\}$ in $V_r$ is adjacent to $4$ vertices $\{u\pm 1,v\}$ and $\{u,v\pm 1\}$. If $n$ is even, all these vertices are in  $V_{r-1}$.
If $n$ is odd, two of them are in $ V_{r-1}$ and the other two in $V_r$.
\end{itemize}
For instance, the quotient graphs of the path-shaped regular partitions of $F_2(C_9)$, $F_2(C_{10})$, and $F_2(C_8)$ are shown in Figures \ref{F_2(C_9)}$(c)$, \ref{F_2(C_{10})}$(c)$, and \ref{F_2(C_8)}$(c)$, respectively.
Then, the quotient matrix $\Q_A=\A(F_2/\pi)$ and quotient Laplacian matrix $\Q_L=\L(F_2/\pi)$ are tridiagonal matrices of the form
\begin{equation}
\label{Qodd}
\Q_A =
\left(
\begin{array}{ccccc}
 0 & 2 &   &   &    \\
 2 & 0 & 2 &  &    \\
    & \ddots & \ddots & \ddots   & \\
     &  & 2 & 0 & 2\\
    &  &  & a & b 
 \end{array}
 \right),\quad \mbox{and}\quad
\Q_L =
\left(
\begin{array}{ccccc}
 2 & -2 &   &   &    \\
 -2 & 4 & -2 &  &    \\
    & \ddots & \ddots & \ddots   & \\
     &  & -2 & 4 & -2\\
    &  &  & -c & c 
 \end{array}
 \right),
 \end{equation}
where $(a,b,c)=(4,0,4)$ if $n$ even, and $(a,b,c)=(2,2,2)$ if $n$ is odd. Then, \eqref{lr-even} and \eqref{lr-odd} correspond to the eigenvalues of the matrix $\Q_L$ for the even and odd cases of $n$, respectively (see Yueh \cite[Th.4]{y05}).
Moreover, in the case of even $n$, the maximum eigenvalue in \eqref{lr-even} is obtained when $r=\nu-1$. Then, from the last comments of Subsection \ref{sec:reg-part}, $\lambda_{\nu-1}$ is the spectral radius of $F_2(C_n)$.
\end{proof}

Next, we show that, depending on the parity of $n$, much more can be stated about the $F_k(C_n)$ spectrum.
Let us consider first the case when $n$ is odd.

\subsection{The case of odd $n$}

As commented in Section \ref{resultats-coneguts}, the case of odd $n$ was studied in Dalf\'o, Fiol, Pavl\'ikov\'a, and \v{S}ir\'an \cite{df22b}, giving the following result.  
\begin{theorem}[\cite{df22b}]
\label{theo:alg-connec}
The 2-token graph $F_2(C_n)$ of the cycle with an odd number $n=2\nu+1$ of vertices is  the lift $G^{\beta}(P^+_{\nu})$ of the base graph the path $P^+_{\nu}$ with vertex set 
$\{u_1,u_2,\ldots,u_{\nu}\}$, a loop at $u_{\nu}$, 
and arcs $a_i = u_iu_{i+1}$ and $a_i^-=u_{i+1}u_i$, for $i = 1,2,\ldots, k-1$.
The voltages on the group $\Z_n$ are as follows:
\begin{align*}
\beta(u_iu_{i+1})&=-1 \qquad \mathrm{for} \  i=1,\ldots,\nu-1,\\
\beta(u_{i+1}u_i)&=+1 \qquad \mathrm{for} \  i=1,\ldots,\nu-1,\\
\beta(u_{\nu}u_{\nu}) &=\pm \nu.
\end{align*}
\end{theorem}

For example, in the case of $n=9$, the 2-token graph $F_2(C_9)$ and its base graph are shown in Figure~\ref{F_2(C_9)}.
Then, the whole spectrum of $F_2(C_n)$ can be obtained from its Laplacian base $\nu\times \nu$ matrix $\B(z)$ (see again \cite{df22b}), with $z=e^{ir\frac{2\pi}{n}}$, or its similar tridiagonal matrix $\B^*(r)$,  for $r=0,1,\ldots,n-1$.
$$
\B(z) =
\left(
{\small{
\begin{array}{cccccc}
2 & -1-z^{-1} & 0 & 0 & \ldots & 0 \\
-1-z & 4 & -1-z^{-1} & 0 & \ldots & 0 \\
0 & -1-z & 4 & -1-z^{-1} & \ddots & 0 \\
0 & 0 & -1-z & \ddots & \ddots  & 0 \\
\vdots  & \vdots & \ddots &\ddots & 4 & -1-z^{-1}\\
0 & 0 & \ldots &0 &-1-z & 4-z^\nu-z^{-\nu} 
\end{array}}}
\right)\cong
$$

{\small{
\begin{align}
 \B^*(r)&=
\left(
\begin{array}{cccccc}
2 & 2\cos(\frac{r\pi}{n}) & 0 & 0& \ldots & 0\\
2\cos(\frac{r\pi}{n}) & 4 & 2\cos(\frac{r\pi}{n})& 0 & \ldots & 0\\
0 & 2\cos(\frac{r\pi}{n}) & 4 & 2\cos(\frac{r\pi}{n})& \ddots& 0\\
0 & 0 & 2\cos(\frac{r\pi}{n}) & \ddots &  \ddots & 0\\
\vdots & \vdots & \ddots & \ddots
& 4 & 2\cos(\frac{r\pi}{n}) \\[.1cm]
0 & 0 & \ldots & 0 & 2\cos(\frac{r\pi}{n}) & 4 + 2(-1)^{r + 1} \cos(\frac{r\pi}{n})
\end{array}
\right).
\label{B-star}
\end{align}
}}
The key idea is that, as shown in the following proposition,  each eigenvector of $\B(z)$, with eigenvalue $\lambda$, gives rise to an eigenvector of $\L$, the Laplacian matrix of $F_2(C_n)$, with the same eigenvalue $\lambda$. Although this follows from results by  
Dalf\'o, Fiol, and \v{S}ir\'a\v{n} \cite{dfs19}, and 
Dalf\'o, Fiol, Pavl\'ikov\'a, and \v{S}ir\'an \cite{df22b}, we give here a direct new proof for completeness. To this end, we label the ${n\choose 2}$ vertices of $F_2(C_n)$, with $n=2\nu+1$, with the pairs $(i,j)=(j+h,j)$ where $j=0,1,\ldots,n-1$ and $h=1,\ldots,\nu$, with arithmetic modulo $n$. With this notation, notice that $h=i-j$ is the distance from $i$ to $j$ in the cycle $C_n$. Besides, every vertex $(i,j)$, with $i=j+h$, of $F_2(C_n)$ corresponds to the vertex $(u_{h},j)$, for $j\in \Z_n$, of the lift $G^{\beta}(P_\nu^{+})$. See Figure \ref{F_2(C_9)}$(a)$-$(b)$ for the case of $F_2(C_9)$ with the simplified notation $ij=(i,j)$.

In what follows, and with a slight abuse of notation, we use $\zeta$ with two meanings. First, $\zeta$ refers to a generic $n$-th root of unity (as in the following proposition). Second, $\zeta^r$ refers to the power $r$ of the $n$-th root of unity so that, in this case, $\zeta=e^{i\frac{2\pi}{n}}$ stands for the first $n$-th root of unity different from $1$ (as in Subsection \ref{sec:sp}).

\begin{proposition}
\label{propo:eigenvec-odd}
Every eigenvalue $\lambda$ of $F_2(C_n)\cong G^{\beta}(P_\nu^{+})$, with $n=2\nu+1$, has an eigenvector $\y\in \Re^{{n\choose 2}}$ with components
\begin{equation}
\label{vector-y}
y_{(i,j)}=f_{i-j}\zeta^j=f_{h}\zeta^j\qquad j=0,\ldots,n-1,\ h=1,\ldots,\nu,
\end{equation}
where $\zeta$ is a given $n$-th root of unity, and $\f= (f_1,\ldots,f_{\nu})$ is a $\lambda$-eigenvector of the matrix $\B(\zeta)$. 
\end{proposition}
\begin{proof}
Let $\L$ be the Laplacian matrix of $F_2(C_n)$. Let $\x$ be an eigenvector of $\L$ with eigenvalue $\lambda$. We show that $\x$ implies the existence of a vector whose entries can be written as claimed. Note first that, since $n$ is odd, all the classes
$V_1,\ldots, V_{\nu}$ of the regular partition of $F_2(C_n)$
have the same number $n$ of vertices. Then, from $\x$, we construct the vectors $\x^0(=\x),\x^1,\ldots,\x^{n-1}$ by shifting in the same way the entries of $\x$ corresponding to each class. More precisely,
\begin{equation}
\label{vector-xa}
x^a_{(i,j)}=x_{(i+a,j+a)}\quad \mbox{ for $a=0,\ldots,n-1$,}
\end{equation}
(all arithmetic understood modulo $n$).
But, in $F_2(C_n)$, the mapping $(i,j)\mapsto (i+a,j+a)$ is an automorphism for every $a=0,\ldots, n-1$. In fact, it is known that
$\Aut F_2(C_n)\cong \Aut C_n$, see Ibarra and Rivera \cite{ir22}. Then, all the vectors $\x^a$ are eigenvectors of $F_2(C_n)$ with eigenvalue $\lambda$. Moreover, from \eqref{vector-xa}, there exists a
matrix $\R$ of size $n\times n$ such that for every $i,j$ (with $i=j+h\ (\mbox{mod } n)$):
\begin{equation}
\label{x-R-x}
\left(
\begin{array}{c}
x^0_{(i+1,j+1)}\\
x^1_{(i+1,j+1)}\\
\vdots \\
x^{n-1}_{(i+1,j+1)}
\end{array}
\right)
= \R
\left(
\begin{array}{c}
x^0_{(i,j)}\\
x^1_{(i,j)}\\
\vdots \\
x^{n-1}_{(i,j)}
\end{array}
\right),
\end{equation}
where 
$$
\R=\circu(0,1,0\ldots,0)=
\left(
\begin{array}{cccccc}
0 & 1 & 0 & \cdots & 0\\
0 & 0 & 1 & \cdots & 0\\
\vdots  & \vdots & \vdots  &\ddots & \vdots \\
 0   & 0 & 0  & \cdots &  1\\
1 & 0 & 0 & \cdots & 0\\
\end{array}
\right).
$$
Then, with $\R^n=\I$, the $n$ eigenvalues of $\R$ are simple  and equal to the $n$-th roots of unity $\zeta^r=e^{i r\frac{2\pi}{n}}$, with $i=\sqrt{-1}$ and $r=0,\ldots,n-1$. Thus,
there exists an inversible (and orthogonal as $\L$ is symmetric) matrix $\Q$ such that $\D=\Q^{-1}\R\Q$, where $\D=\diag(1,\zeta,\zeta^2,\ldots,\zeta^{n-1})$. Let $\X$ be the matrix with rows $\x^0,\x^1,\ldots,\x^{n-1}$. From $\X\L=\lambda\X$, we have $\Q^{-1}\X\L=\lambda\Q^{-1}\X$ and, hence, the vectors $\y^0,\y^1,\ldots,\y^{n-1}$ with components
\begin{equation}
\label{y-Q-x}
\left(
\begin{array}{c}
y^0_{(i,j)}\\
y^1_{(i,j)}\\
\vdots \\
y^{n-1}_{(i,j)}
\end{array}
\right)
= \Q^{-1}
\left(
\begin{array}{c}
x^0_{(i,j)}\\
x^1_{(i,j)}\\
\vdots \\
x^{n-1}_{(i,j)}
\end{array}
\right),
\end{equation}
or rows of the matrix $\Y=\Q^{-1}\X$,
are also $\lambda$-eigenvectors of $F_2(G)$, provided that they are different from $\vec0$. (Notice that the number of such eigenvectors, that is, $\rank\Y$, is at most $m(\lambda)$, the multiplicity of $\lambda$.) Moreover, from \eqref{y-Q-x} and \eqref{x-R-x} with $\R=\Q\D\Q^{-1}$, we get 
\begin{equation}
\label{y-R-y}
\left(
\begin{array}{c}
y^0_{(i+1,j+1)}\\
y^1_{(i+1,j+1)}\\
\vdots \\
y^{n-1}_{(i+1,j+1)}
\end{array}
\right)
= \D
\left(
\begin{array}{c}
y^0_{(i,j)}\\
y^1_{(i,j)}\\
\vdots \\
y^{n-1}_{(i,j)}
\end{array}
\right).
\end{equation}
Now, since $\x\neq \vec0$, there is at least a non-zero vector $\y^a$ satisfying
$y^a_{(i+1,j+1)}=\zeta^a y^a_{(i,j)}$ or, iterating, $y^a_{(i+b,j+b)}=\zeta^{ab} y^a_{(i,j)}$.
Therefore, letting $b=-j$,
$$
y^a_{(i,j)}= y^a_{(i-j,0)}\zeta^{aj}.
$$
Then, the result in \eqref{vector-y} follows by taking $\y=\y^a$, $\zeta=\zeta^a$, and $f_h=y^a_{(h,0)}$.
To show that $\f(\zeta)$ is an eigenvector of $\B(z)$ with $z=\zeta$, we distinguish three cases:
 \begin{itemize}
\item[$(i)$] $h=1$: The vertex $(j+1,j)$, with vector entry $y_{(j+1,j)}=f_1\zeta^j$, is adjacent to both vertices $(j+2,j)$ and $(j+1,j-1)$, with respective vector entries $f_2\zeta^j$ and $f_2\zeta^{j-1}$. Then, the first entry of $\y\L=\lambda\y$ is $2f_1\zeta^j-f_2\zeta^j-f_2\zeta^{j-1}=\lambda f_1\zeta^j$. Hence, $2f_1-f_2-f_2\zeta^{-1}=\lambda f_1$, which corresponds to the first entries of $\B(\zeta)\f=\lambda\f$, as it should.
\item[$(ii)$] $2\le h\le \nu-1$: The vertex $(j+h,j)$, with vector entry $y_{(j+h,j)}=f_h\zeta^j$, is adjacent to the four vertices $(j+h+1,j)$, $(j+h,j+1)$, $(j+h-1,j)$, and $(j+h,j-1)$, with respective vector entries $f_{h+1}\zeta^j$, $f_{h-1}\zeta^{j+1}$, $f_{h-1}\zeta^{j}$, and $f_{h+1}\zeta^{j-1}$. Then, the $h$-th entry of $\y\L=\lambda\y$ is 
 $$
 4f_h\zeta^j-f_{h+1}\zeta^j-f_{h-1}\zeta^{j+1}-f_{h-1}\zeta^{j}-f_{h+1}\zeta^{j-1}=\lambda f_h\zeta^j.
 $$ 
Thus,
$$
4f_h-f_{h+1}-f_{h-1}\zeta-f_{h-1}-f_{h+1}\zeta^{-1}=\lambda  f_h,
$$ 
which corresponds to the $h$-th entry of $\B(\zeta)\f=\lambda\f$.

\item[$(iii)$] $h=\nu$: 
The vertex $(j+\nu,j)$, with vector entry $y_{(j+\nu,j)}=f_{\nu}\zeta^j$,  is adjacent (according to the used notation) to the four vertices $(j,j+\nu+1)$, $(j+\nu-1,j)$, $(j+\nu,j+1)$, and $(j-1,j+\nu)$, with respective vector entries $f_{-(\nu+1)}\zeta^{j+\nu+1}=f_{\nu}\zeta^{j+\nu+1}$ (the subscripts of $f$ are modulo $n$), $f_{\nu-1}\zeta^{j}$, $f_{\nu-1}\zeta^{j+1}$, and $f_{-(1+\nu)}\zeta^{j+\nu}=f_{\nu}\zeta^{j+\nu}$. Then, the $\nu$-th entry of $\y\L=\lambda\y$ is 
 $$
 4f_{\nu}\zeta^j-f_{\nu}\zeta^{j+\nu+1}-f_{\nu-1}\zeta^{j}-f_{\nu-1}\zeta^{j+1}-f_{\nu}\zeta^{j+\nu}=\lambda f_{\nu}\zeta^j.
 $$ 
Thus,
$$
4f_{\nu}-f_{\nu}\zeta^{\nu+1}-f_{\nu-1}-f_{\nu-1}\zeta-f_{\nu}\zeta^{\nu}=\lambda  f_{\nu},
$$ 
which corresponds to the $\nu$-th entry of $\B(\zeta)\f=\lambda\f$.

\end{itemize}
This completes the proof.
\end{proof}

In the case of $F_2(C_9)$, the obtained eigenvalues are shown in  Table \ref{taula:C9}.


We focus on the matrix $\B^*(r)$ in the following result.
\begin{proposition}
\label{propo:B*(r)}
Given $r=0,1,\ldots,2\nu$, let $\lambda_{r,1}\le\lambda_{r,2}\le\cdots \le\lambda_{r,\nu}$ be the eigenvalues of the matrix $\B^*(r)$. Then, the following holds:
\begin{itemize}
\item[$(i)$]
The eigenvalues of $\B^*(0)$ are
\begin{equation*}
\lambda_{0,s} =8\cos^2 \left(\frac{s\pi}{2\nu}\right)\quad \mbox{for $s=1,2,\ldots, \nu$}. \label{lj}
\end{equation*}
Then, $\lambda_{0,\nu}=0$. Moreover, the smallest nonzero eigenvalue is obtained when $s=\nu-1$ and  satisfies 
$\lambda_{0,\nu-1}>\alpha(C_{2\nu+1})$.
\item[$(ii)$] 
For each $r=1,\ldots,2\nu$, the smallest 
eigenvalue of $\B^*(r)$ satisfies
$$
\lambda_{r,1}\ge 4\sin^2\left(\frac{r\pi}{2(2\nu+1)}\right).
$$
\item[$(iii)$]
The matrix $\B(\zeta^r)$, with $\zeta=e^{i\frac{2\pi}{2\nu+1}}$,  has exactly one eigenvalue of $C_{2\nu+1}$, which is $\lambda_r=2\left(1-\cos\left(\frac{r2\pi}{2\nu+1}\right)\right)$.
Besides, the eigenvalues of $\B^*(0)$ have multiplicity one, whereas  the eigenvalues of $\B^*(r)$, with $r=1,\ldots, 2\nu$, have multiplicity two. 

\end{itemize}
\end{proposition}
\begin{proof}
$(i)$ Notice that, when $r=0$, $\B^*(0)$ is the tridiagonal matrix $\Q_L$ in \eqref{Qodd}, with eigenvalues already given in Proposition \ref{propo:path-shaped}.
Moreover, the function 
$$
\phi(\nu)=\frac{\alpha(C_{2\nu+1})}{\lambda_{0,\nu-1}}=
\frac{\sin^2\left(\frac{\pi}{2\nu+1}\right)}{2\cos^2\left(\frac{(\nu-1)\pi}{2\nu}\right)}
$$
satisfies $\phi(\nu)< 1/2$ for  $\nu>0$.
\\
$(ii)$
From the matrix $\B^*(r)$, we have three different Gershgorin circles determining three intervals $I_i(r)$, for $i=1,2,3$, in the real line. Thus, all eigenvalues of $\B^*(r)$ are within $I_1(r)\cup I_2(r)\cup I_3(r)$.   
The left endpoints of these intervals are $\ell_1(r)=4\sin^2(\frac{r\pi}{2(2\nu+1)})$, $\ell_2(r)=8\sin^2(\frac{r\pi}{4\nu+2})$, and
$\ell_3(r)=4$ for $r$ odd and 
$\ell_3(r)=8\sin^2(\frac{r\pi}{4\nu+2})$ for $r$ even.
Then, 
\begin{align*}
\ell_i(0) &=0\quad \mbox{for $i=1,2,3$,}\\
\ell_1(1) &=4\sin^2\left(\frac{\pi}{2(2\nu+1)}\right)=\alpha(P_{2\nu+1}),\\ 
\ell_1(2) &=4\sin^2\left(\frac{\pi}{2\nu+1}\right)=\alpha(C_{2\nu+1}).
\end{align*}
(See Table \ref{taula:limits-C9} for the case $\nu=4$, corresponding to the matrices $\B^*(r)$ of $F_2(C_9)$.)
Now, a simple analysis shows that the values of $\ell_i(r)$ are increasing when $r=0,1,\ldots, 2\nu$ when $i=1,2$, and $r=0,2,\ldots, 2\nu$ when $i=3$. Thus, the result follows since
$\ell_1(r)<\min\{\ell_2(r),\ell_3(r)\}$ for $r\neq 0$.

$(iii)$ We prove that, for every $z=e^{ir\frac{2\pi}{n}}$, for $r=0,\ldots,n-1$, the matrix $\B(z)\cong \B^*(r)$ has exactly one eigenvalue of $C_n$.
From Proposition \ref{propo:eigenvec-odd}, we know that each of the eigenvectors $\y$ of $F_2(C_{2\nu+1})$ 
has entries 
$y_{(i,j)}=f_{i-j}\zeta^j=f_{h}\zeta^j$ for $i=0,\ldots,n-1$, $h=1,\ldots,\nu$, and $\f=(f_1,\ldots,f_{\nu})$ an eigenvector of $\B(\zeta)$.
 Let $\B$ be the $(n,2)$-binomial matrix. Recall that, from Lemma \ref{coro:LkL1}$(ii)$,  if $\y$ is a $\lambda$-eigenvector of $F_2(G)$ and $\x=\B^{\top}\y\neq \vec0$, then $\x$ is a $\lambda$-eigenvector of $G$.
In our case, notice that, for $j=0,\ldots,n-1$, the $j$-th entry of  the vector $\x=\B^{\top}\y$
is obtained by adding all the $n-1$ entries of the vector $\y$ with labels having
a common $j$, that is, corresponding to the vertices 
$$
(j+1,j),(j+2,j),\ldots,(j+\nu,j),(j,j-\nu),(j,j-\nu+1),\ldots,(j,j-1).
$$
Then, the $j$-th entry of the vector $\x$ turns out to be
\begin{equation}
\label{x-from-f}
f_1(\zeta^j+\zeta^{j-1})+f_2(\zeta^j+\zeta^{j-2})+\cdots + 
f_{\nu}(\zeta^j+\zeta^{j-\nu}),\qquad j=0,\ldots,n-1.
\end{equation}
Let $\F$ be the $\nu\times \nu$  matrix whose columns $\f^1,\ldots,\f^{\nu}$ are the eigenvectors of $\B(\zeta)$, and $\ZZ$ the $n\times \nu$ matrix with $j$-th row $(\zeta^j+\zeta^{j-1}), (\zeta^j+\zeta^{j-2}),\ldots,
(\zeta^j+\zeta^{j-\nu})$.
Then, in matrix form, \eqref{x-from-f} is $\X=\ZZ\F$, where $\X$ is the $n\times \nu$ matrix with columns being the putative eigenvectors of $C_n$. But $\F$ has full rank, whereas $\ZZ$
has rank 1 (every row is a multiple of the first one). Consequently, $\rank \X=1$ and, hence, exactly one $\lambda$-eigenvector of $\B(\zeta)$ gives a $\lambda$-eigenvector of $C_n$.
The statement about the multiplicities follows from the fact that if $\zeta\neq 1$, the spectra of $\B(\zeta)$ and $\B(\zeta^{-1})=\B(\overline{\zeta})$ coincide, where $\overline{\zeta}$ is the conjugate of $\zeta$.
\end{proof}

\begin{table}[t!]
\begin{center}
\begin{tabular}{|c|cccc| }
\hline
$\zeta=e^{i\frac{2\pi}{9}}$, $z=\zeta^r$ & $\lambda_{r,1}$  & $\lambda_{r,2}$  & $\lambda_{r,3}$  & $\lambda_{r,4}$   \\
\hline\hline
$\spec(\B(\zeta^0))$ & \bf 0 & 1.171572876 & 4 &  6.828427124   \\
\hline
$\spec(\B(\zeta^1))=\spec(\B(\zeta^8))$ & \bf 0.4679111136 & 2.52079560 & 5.420264509 &  7.470414013   \\
\hline
$\spec(\B(\zeta^2))=\spec(\B(\zeta^7))$ & 0.783324839 & \bf 1.65270363  &  3.895673125 & 6.136209510  \\
\hline
$\spec(\B(\zeta^3))=\spec(\B(\zeta^6))$ & 1.50913638   & \bf 3  &  4.656620432 &  5.834243185  \\
\hline
$\spec(\B(\zeta^4))=\spec(\B(\zeta^5))$ & 1.939683655   & 3.382489411  & \bf  3.87938479 &  4.451145779  \\
\hline
\end{tabular}
\end{center}
\caption{All the eigenvalues of the matrices $\B(\zeta^r)$, which yield the eigenvalues of the 2-token graph $F_2(C_9)$. The values in boldface correspond to the eigenvalues of $C_9$.}
\label{taula:C9}
\end{table}

\begin{table}[t!]
\begin{center}
\begin{tabular}{|c|cccc| }
\hline
  & $r=0$  & $r=1$  & $r=2$ & $r=3$ \\
\hline\hline
$\ell_1(r)$ & 0 & \bf 0.1206147584 & \bf 0.4679111138 & 1  \\
\hline
$\ell_2(r)$ & 0 & 0.24122951686 & 0.93582222752 & 2 \\
\hline
$\ell_3(r)$ & 0 & 4  &  0.93582222752 & 4\\
\hline
\end{tabular}
\end{center}
\caption{Left endpoints of the Gershgorin circles of the matrices $\B^*(r)$.
The values in boldface correspond to the algebraic connectivities of $P_9$ and $C_9$.}
\label{taula:limits-C9}
\end{table}

Moreover, results from Table \ref{taula:C9} suggest that the minimum and maximum eigenvalues of the matrix $\B^*(1)$
correspond to the algebraic connectivity $\alpha(F_2(C_{2\nu+1}))=\alpha(C_{2\nu+1})$, and
spectral (Laplacian) radius $\rho(F_2(C_{2\nu+1}))$, respectively.


\subsection{The case of even $n$ and odd $n/2$}

In the case of even $n=4r+2$ (so that $n/2$ is odd), the $2$-token $F_2(C_n)$ can also be seen as a lift graph, as shown in the following result. See an example of this kind of token graph in Figure \ref{F_2(C_{10})}.


\begin{figure}[t]
\begin{center}
\includegraphics[width=14cm]{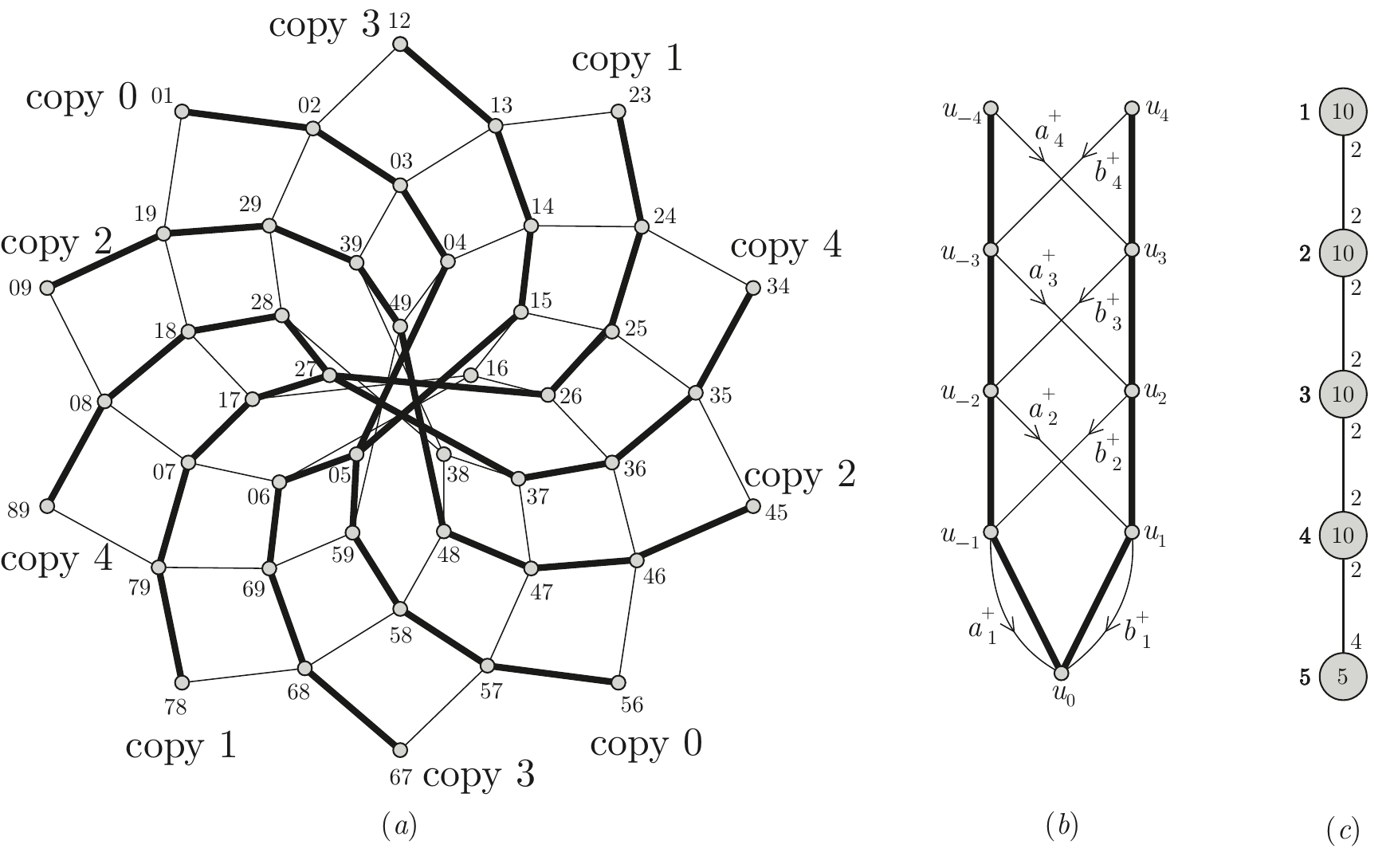}
\caption{$(a)$ The token graph $F_2(C_{10})$ of the cycle graph $C_{10}$, with the different copies of the U-shaped regular partition, here drawn as paths. $(b)$ Its base digraph with voltages in $\mathbb{Z}_{10}$, which gives rise to the U-shaped regular partition. The thick edges represent the path $P_9$ obtained with this partition. $(c)$ The quotient graph of the path-shaped regular partition. In boldface, there is the numbering of the vertex classes.}
\label{F_2(C_{10})}
\end{center}
\end{figure}

Given an integer $r$, let us consider the path graph $G=P_{4r+1}$ with vertices 
$$
u_{-2r},u_{-2r+1},\ldots, u_{-1},u_0,u_1,\ldots,u_{2r-1},u_{2r},
$$
(with its corresponding edges) and additional arcs 
$$
\begin{array}{lll}
a_i^{+}=u_{-i}u_{i-1}, & a_i^{-}=u_{i-1}u_{-i}, & \mbox{for } i=0,1\ldots,2r,\\
b_i^{+}=u_{i}u_{-i+1}, & b_i^{-}=u_{-i+1}u_{i}, & \mbox{for } i=0,1,\ldots,2r.
\end{array}
$$


Let $\beta$ be the voltage assignment  on $G$ in the cyclic group $\Z_{2r+1}$ given by 
\begin{equation}
\label{eq:voltatgesC10}
    \begin{array}{c}
\beta(a_i^+)=\beta(b_i^+)=+r,\\
\beta(a_i^-)=\beta(b_i^-)=-r.
\end{array}
\end{equation}
Figure \ref{F_2(C_{10})}$(b)$ shows the base graph $G$ for $r=2$. Now, we have the following result.

\begin{lemma}
Given $G=P_{4r+1}$ with the voltage assignment \eqref{eq:voltatgesC10} on $\Z_{2r+1}$, the 2-token graph of the cycle $C_{n}$ with $n=4r+2$ is the lift graph $G^{\beta}$. That is,
$$
F_2(C_{4r+2})\cong G^{\beta}.
$$
\end{lemma}

\begin{proof}
The vertex set $V(G^{\beta})$ of the lift $G^{\beta}$ has elements labeled with the pairs $(u_i,g)$
for $i=-2r,\ldots,2r$ and $g\in \Z_{2r+1}$. Thus $|V(G^{\beta})|=(4r+1)(2r+1)={4r+2\choose 2}$, which corresponds to the number of vertices of $F_2(C_{4r+2})$. Indeed, such vertices correspond to 2-subsets $\{i,j\}$ of the set $\{1,2,\ldots,4r+2\}$.
Thus, we have to show a 1-to-1 mapping between $V(G^{\beta})$ and $V(F_2(C_{4r+2}))$ that must be consequent with the adjacencies of both graphs. Such a mapping is shown in Table \ref{taula:token-lifts}, from where it is easily checked that the adjacency conditions are fulfilled.
Let us take an example. The vertex $(u_{-2r+1},1)\equiv \{2,4\}$ (written as 24 in Figure \ref{F_2(C_{10})}) of $G^{\beta}$
is adjacent to:
\begin{itemize}
    \item 
    The vertices $(u_{-2r},1)$ and $(u_{-2r+2},1)$ of the same `copy'.
    \item 
    The vertex  $(u_{2r-2},r+1)$ by the arc $a_{2r-1}^+$ with voltage $+r$.
    \item
    The vertex 
    $(u_{2r},r+2)$ by the arc $b^{-}_{2r}$ with voltage $-r$.
\end{itemize}
Then, looking again at Table \ref{taula:token-lifts}, we find the following equivalences:
\begin{itemize}
    \item 
   $(u_{-2r},1)\equiv \{2,3\}$ and $(u_{-2r+2},1)\equiv \{2,5\}$,
    \item 
    $(u_{2r-2},r+1)\equiv \{1,4\}$ and $(u_{2r},r+2)\equiv \{3,4\}$,
\end{itemize}
which correspond to the vertices adjacent to $\{2,4\}$ in $F_2(C_{4r+2}$).

\begin{table}[t!]
\begin{center}
\scriptsize
\begin{tabular}{|c|c@{}c@{}c@{}c@{\hspace{3pt}}c@{\hspace{2pt}}c@{\hspace{2pt}}| }
\hline
 \backslashbox{Vertex}{Copy $g$} & $0$ & $r+1$ & $1$ & $r+2$ & $\cdots$ & $r$\\
\hline
$u_{-2r}$ & $\{0,1\}$ & $\{2r+2,2r+3\}$ & $\mathbf{\{2,3\}}$ & 
$\{2r+4,2r+5\}$ & $\cdots$ & 
$\{2r,2r+1\}$  \\
$u_{-2r+1}$ & $\{0,2\}$ & $\{2r+2,2r+4\}$ & $\mathbf{\{2,4\}}$ & 
$\{2r+4,2r+6\}$ & $\cdots$ & 
$\{2r,2r+2\}$ \\
$u_{-2r+2}$ & $\{0,3\}$ & $\{2r+2,2r+5\}$ & $\mathbf{\{2,5\}}$ & 
$\{2r+4,2r+7\}$ & $\cdots$ & 
$\{2r,2r+3\}$ \\
$\vdots$ & $\vdots$ & $\vdots$ & $\vdots$ & $\vdots$ & &$\vdots$ \\
$u_{-1}$ & $\{0,2r\}$ & $\{2r+2,0\}$ & $\{2,2r+2\}$ & $\{2r+4,2\}$ & $\cdots$ & 
$\{2r,4r\}$ \\
$u_{0}$ &$\{0,2r+1\}$ & $\{1,2r+2\}$ & $\{2,2r+3\}$ & $\{3,2r+4\}$ & $\cdots$ & 
$\{2r,4r+1\}$ \\
$u_{1}$ & $\{2r+1,4r+1\}$ & $\{1,2r+1\}$ & $\{2r+3,1\}$ & 
$\{3,2r+3\}$ & $\cdots$ & $\{4r+1,2r-1\}$ \\
$\vdots$ & $\vdots$ & $\vdots$ & $\vdots$ & $\vdots$ & & $\vdots$  \\
$u_{2r-2}$ & $\{2r+1,2r+4\}$ & $\mathbf{\{1,4\}}$ & $\{2r+3,2r+6\}$ & 
\{3,6\} & $\cdots$ & $\{4r+1,2\}$ \\
$u_{2r-1}$ & $\{2r+1,2r+3\}$ & $\{1,3\}$ & $\{2r+3,2r+5\}$ & 
\{3,5\} & $\cdots$ & $\{4r+1,1\}$ \\
$u_{2r}$ & $\{2r+1,2r+2\}$ & $\{1,2\}$ & $\{2r+3,2r+4\}$ & 
$\mathbf{\{3,4\}}$ & $\cdots$ & $\{4r+1,0\}$\\
\hline
\end{tabular}
\end{center}
\caption{All the vertices of the 2-token graph $F_2(C_n)$ for even $n$ and odd $n/2$. The vertex \{2,4\} and its adjacent vertices are in boldface, following the example.}
\label{taula:token-lifts}
\end{table}
\end{proof}

As a consequence, for each $z=\zeta ^{\ell}$, where $\zeta= e^{i\frac{2\pi}{2r+1}}$, with $\ell=0,1,\ldots,2r$,  an irreducible representation of  the Laplacian base matrix of $G^{\beta}\cong F_2(C_{2r+2})$ is the matrix $\B(z)$ 
as shown next.

Note that matrix $\B(z)$ is tridiagonal with respect to the main and the secondary diagonals.

$$
     \B(z)=
     \scriptsize
       \left(
       \arraycolsep=1pt
\begin{array}{ccccccccccccc}
2      & -1    & 0     & \cdots& 0     & 0       & 0      & 0       & 0     & \cdots& 0     & -z^r  & 0 \\
-1     & 4     & -1    & \cdots& 0     & 0       & 0      & 0       & 0     & \cdots&-z^{r} & 0     & -z^{-r}\\
0      &-1     & 4     & \cdots& 0     & 0       & 0      & 0       & 0     &\cdots & 0     &-z^{-r}& 0 \\
\vdots & \vdots& \vdots& \ddots& \vdots& \vdots  & \vdots & \vdots  & \vdots&\iddots&\vdots & \vdots& \vdots \\
0      & 0     & 0     & \cdots& 4     & -1      & 0      &-z^{r}   & 0     &\cdots & 0     & 0     & 0 \\
0      & 0     & 0     & \cdots& -1    & 4       &-1-z^{r}& 0       &-z^{-r}& \cdots& 0     & 0     & 0 \\
0      & 0     & 0     & \cdots& 0     &-1-z^{-r}& 4      &-1-z^{-r}& 0     & \cdots& 0     & 0     & 0\\
0      & 0     & 0     & \cdots&-z^{-r}& 0       &-1-z^{r}& 4       &-1     & \cdots& 0     & 0     & 0 \\
0      & 0     & 0     & \cdots& 0     & -z^{r}  & 0      & -1      & 4     & \cdots& 0     &  0    & 0  \\
\vdots & \vdots&\vdots &\iddots& \vdots& \vdots  & \vdots & \vdots  & \vdots& \ddots& \vdots& \vdots& \vdots\\
0      &-z^{-r}& 0     &\cdots & 0     & 0       & 0      & 0       & 0     & \cdots& 4     & -1    & 0\\
-z^{-r}& 0     &-z^{r} &\cdots & 0     & 0       & 0      & 0       & 0     & \cdots& -1    & 4     & -1 \\
0      &-z^r   & 0     & \cdots& 0     & 0       & 0      & 0       & 0     & \cdots& 0     & -1    & 2 \\
\end{array}
\right).  
$$

For instance, in the case of  $F_2(C_{6})$ ($r=1$), we have
$$
\B(z)=
\left(
\begin{array}{ccccc}
2 & -1 & 0 & -z & 0 \\
-1 & 4 & -1-z & 0 & -z^{-1} \\
0 & -1-z^{-1} & 4 & -1-z^{-1}& 0 \\
-z^{-1} &  0 & -1-z & 4 & -1\\
0 & -z & 0 & -1 & 2
\end{array}
\right).
$$
In Table \ref{table5}, we show the different eigenvalues of $F_2(C_{6})$, obtained as the eigenvalues of each $\B(z)$ for $z=\zeta^{\ell}$ with $\ell =0,1,2$.

\begin{table}[ht]
\begin{center}
\begin{tabular}{|c|ccccc| }
\hline
$\zeta=e^{i\frac{2\pi}{3}}$, $z=\zeta^r$ & $\lambda_{r,1}$  & $\lambda_{r,2}$  & $\lambda_{r,3}$  & $\lambda_{r,4}$ & $\lambda_{r,5}$   \\
\hline\hline
$\spec(\B(\zeta^0))$ & \bf 0 & 2 & $5-\sqrt{5}\approx 2.764$ &  \bf 4 & $5+\sqrt{5}\approx 7.236$  \\
\hline
$\spec(\B(\zeta^1))$ & \bf 1 & $\frac{1}{2}(7-\sqrt{17})\approx 1.438$ & \bf 3 &  5 & $\frac{1}{2}(7+\sqrt{17})\approx 5.561$   \\
\hline
$\spec(\B(\zeta^2))$ &   \bf 1 & $\frac{1}{2}(7-\sqrt{17})\approx 1.438$ & \bf 3 &  5 & $\frac{1}{2}(7+\sqrt{17})\approx 5.561$   \\
\hline
\end{tabular}
\end{center}
\caption{All the eigenvalues of matrices $\B(\zeta^r)$, which yield the eigenvalues the 2-token graph  $F_2(C_6)$. The values in boldface correspond to the eigenvalues of $C_6$.}
\label{table5}
\end{table}

\subsection{The case of even $n$ and $n/2$}
\begin{figure}[t]
\begin{center}
\includegraphics[width=14cm]{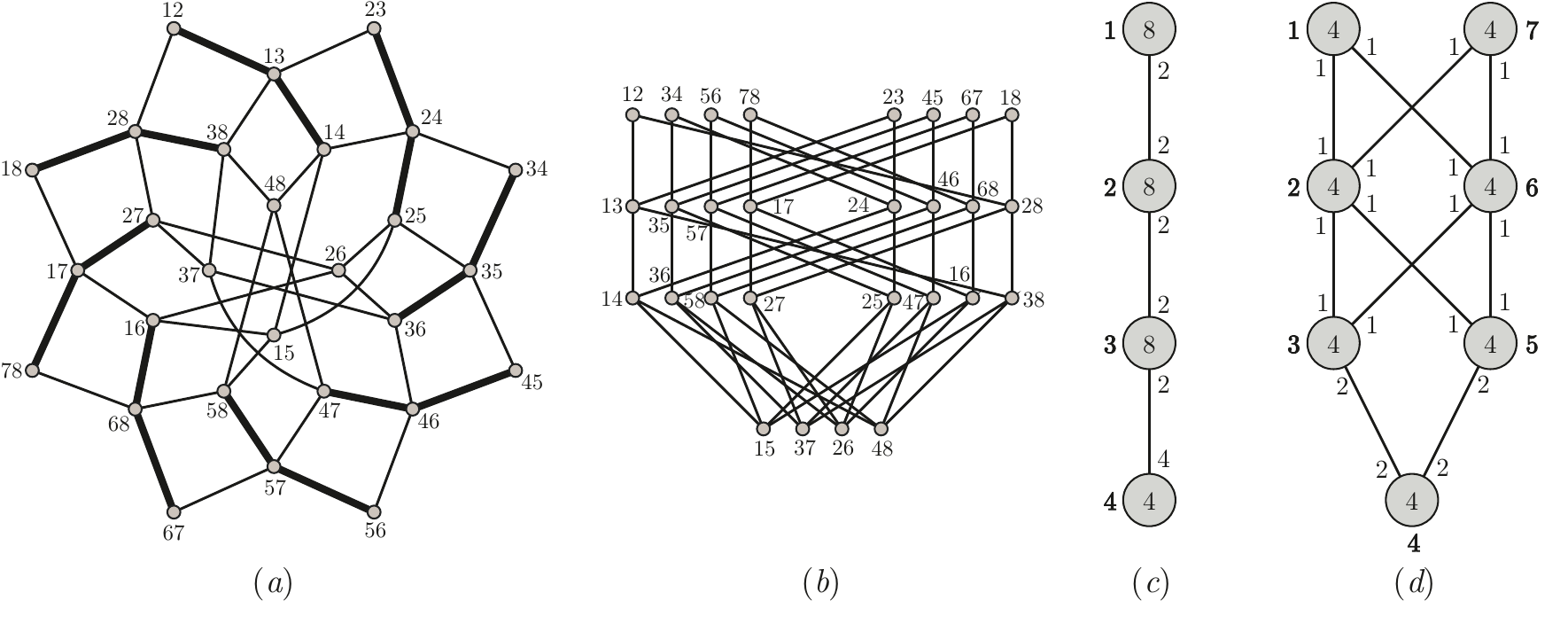}
\caption{$(a)$ The $2$-token graph $F_2(C_8)$ of the cycle graph on 8 vertices. $(b)$ Another view of $F_2(C_8)$. $(c)$ The quotient graph from the path-shaped regular partition. $(d)$ The quotient graph of the U-shaped regular partition obtained from $(b)$. In boldface, there is the numbering of the vertex classes.}
\label{F_2(C_8)}
\end{center}
\end{figure}

When we consider the case of cycles $C_n$ with $n$ and $n/2$ even, it is not very useful to represent $F_2(C_n)$ as a lift graph to compute the whole spectrum. The reason is that the base graph has too many vertices with respect to the original graph. Alternatively, besides the spectrum of $C_n$, we can easily find another part of the spectrum by means of regular partitions. In fact, we can use the regular path-partition and U-partition with the same structure as in the previous subsection.   As an example, the 2-token graph of  $C_8$ is shown in Figure \ref{F_2(C_8)}$(a)$ and $(b)$, together with its regular partitions $(c)$ (the path $P_{n/2}$) and $(d)$ (the U-shaped graph). Compare such partitions with those in  Figure \ref{F_2(C_{10})}$(b)$ and $(c)$.
Then, the quotient $7\times 7$ Laplacian matrix of the regular U-shaped partition $\tau$ of $F_2=F_2(C_8)$ and its spectrum are
\begin{align*}
\Q(F_2/\tau) &=
\left(
\begin{array}{ccccccc}
2 & -1 & 0 & 0 & 0 & -1& 0\\
-1 & 4 & -1 & 0 & -1 &  0& -1\\
0 & -1 & 4 & -2 & 0 & -1 &  0\\
0 & 0 & -2 & 4 & -2 & 0 & 0\\
0 & -1 & 0 & -2 & 4 & -1& 0\\
-1 & 0 & -1 & 0 & -1 & 4& -1\\
0 & -1 & 0 & 0 & 0 & -1 & 2
\end{array}
\right),\\
\spec \Q(F_2/\tau)  &=\{0,1.5060,2,4^2,4.8900,7.6038\} \ \subset \spec F_2(C_{8}).
\end{align*}

With respect to the path-shaped partition $\pi$, Proposition \ref{propo:path-shaped} gives that the quotient Laplacian matrix and its spectrum are
\begin{align*}
\Q(F_2/\pi) &=
\left(
\begin{array}{cccc}
2 & -2 & 0 & 0 \\
-2 & 4 & -2 & 0 \\
0 & -2 & 4 & -2\\
0 & 0 & -4 & 4\\
\end{array}
\right),\\
\spec \Q(F_2/\pi)  &=\{0,1.5060,4.8900,7.6038\} \ \subset \spec \Q(F_2/\tau) \ \subset \spec F_2(C_8).
\end{align*}
Notice that the inclusion $\spec \Q(F_2/\pi) \subset \spec \Q(F_2/\tau)$ is due to the fact that $\tau$ can be seen as a regular partition of the quotient graph $F_2/\pi$.
Moreover, according to the same proposition,  the largest eigenvalue $7.6038$ of $\Q(F_2/\pi)$ (or $\Q(F_2/\tau)$) is the spectral radius $\rho(F_2(C_8))$.

Compare these results with the whole spectrum of $F_2(C_8)$, which is
\begin{align}
\spec F_2(C_8) =\{&0,0.5857^2,0.9486^2,1.5060,1.7117^2,2^3,3.1259^2,3.4142^2, \nonumber\\
&4^3,4.5173^2,4,8740^2,4.8900,6.2882^2,6.5340^2,7.6038\}.
\label{spF2(C8)}
\end{align}

\subsection{The new method of over-lifts}
In this subsection, we use a new method called {\em over-lifts}, which allows us
to unify the cases of cycles with even $n$ and compute the whole spectrum of $F_2(C_n)$. 
(For instance, as a result of such a method, all the eigenvalues in  \eqref{spF2(C8)} are shown in Table \ref{taula:C8}.)
This is accomplished by means of a new polynomial matrix $\B(z)$ that does {\bf not} correspond to the base graph of a lift. By its characteristics, we say that $\B(z)$ is associated with an over-lift. The basic difference is that such a matrix has dimension $\nu\times\nu$ (recall that $n=2\nu$), and there are $n$ possible values for $z$ ($n$-th roots of unity).
Thus, the total number of eigenvalues obtained is $\nu n$. However, $F_2(C_n)$ has ${n\choose 2}=\nu(n-1)$ vertices, which is the number of eigenvalues of $\L$. 
We will see that, in fact, the $\nu$ `extra' eigenvalues provided by $\B(z)$ are all equal to $4$.

\begin{theorem}
Let $\L$ be the Laplacian matrix of $F_2(C_n)$, with $n=2\nu$.
Let $\Lambda$ be the multiset with elements $4,4,\stackrel{(\nu)}{\ldots},4$. Then,
the spectrum of $\L$ can be obtained from the spectrum of the $\nu\times\nu$ matrix $\B(z)$ or, equivalently, from the spectrum of its similar matrix $\B^*(r)$:
$$
\B(z) =
\left(
\small{
\begin{array}{cccccc}
2 & -1-z^{-1} & 0 & 0 & \ldots & 0 \\
-1-z & 4 & -1-z^{-1} & 0 & \ldots & 0 \\
0 & -1-z & 4 & -1-z^{-1} & \ddots & 0 \\
0 & 0 & -1-z & \ddots & \ddots  & 0 \\
\vdots  & \vdots & \ddots &\ddots & 4 & -1-z^{-1}\\
0 & 0 & \ldots &0 &-1-z-z^{\nu}-z^{\nu+1} & 4 
\end{array}
}
\right)
$$
for $z=\zeta^r=e^{ir\frac{2\pi}{n}}$, $r=0,1,\ldots,n-1$, and 

\small{
\begin{align}
 \B^*(r)&=
\left(
\begin{array}{cccccc}
2 & 2\cos(\frac{r\pi}{n}) & 0 & 0& \ldots & 0\\
2\cos(\frac{r\pi}{n}) & 4 & 2\cos(\frac{r\pi}{n})& 0 & \ldots & 0\\
0 & 2\cos(\frac{r\pi}{n}) & 4 & 2\cos(\frac{r\pi}{n})& \ddots& 0\\
0 & 0 & 2\cos(\frac{r\pi}{n}) & \ddots &  \ddots & 0\\
\vdots & \vdots & \ddots & \ddots
& 4 & 2\cos(\frac{r\pi}{n}) \\[.1cm]
0 & 0 & \ldots & 0 & 2\cos(\frac{r\pi}{n})+2\cos(\frac{r(n-1)\pi}{n}) & 4 
\end{array}
\right)
\label{B-star-even}
\end{align}
for $r=0,1,\ldots,n-1$. 
Formally,  
\begin{equation}
\label{spB(z)-spL}
\bigcup_{z\in R(n)}\spec \B(z)=\bigcup_{r=0}^{n-1}\spec \B^*(r)= \Lambda \cup\spec \L, 
\end{equation}
where $R(n)$ denotes the set of the $n$-th roots of unity,
}
\end{theorem}
\begin{proof}
For convenience, we will use indistinctly $\B(z)$ or $\B^*(r)$
because 
of $z=\zeta^r$. 
First, let us prove that, as in the case of odd $n$, most of the eigenvectors 
$\f=(f_1,f_2,\ldots,f_{\nu})$ of $\B(z)$, with eigenvalue $\lambda$, give rise to an eigenvector of  $\L$, with the same eigenvalue $\lambda$.
Indeed, with the same notation for the vertices $(i,j)=(i+h,j)$ of $F_2(C_n)$ as in Proposition \ref{propo:eigenvec-odd}, let us consider the vector $\x\in \Re^{{n\choose 2}}$ with components
\begin{equation}
x_{(i,j)}=f_{h}\zeta^j,\qquad i=0,\ldots,n-1, \ h=1,\ldots,\nu,
\end{equation}
where $\zeta$ is a given $n$-th root of unity.
Now, to show that, under some conditions, $\f$ is an eigenvector of $\B(z)$ with $z=\zeta$, we distinguish four cases: $(i)$ $h=1$; $(ii)$ $2\le h\le \nu-2$; $(iii)$ $h=\nu-1$; and $(iv)$ $h=\nu$. Since the cases $(i)$ and $(ii)$ are proved as in Proposition \ref{propo:eigenvec-odd}, we only consider $(iii)$ and $(iv)$.
Since case $(iii)$ is the most involved, we begin with $(iv)$.
 \begin{itemize}
 \item[$(iv)$] $h=\nu$: 
The vertex $(j+\nu,j)=(j,j+\nu)$, with vector entry $x_{(j+\nu,j)}=f_{\nu}\zeta^j$, is adjacent to the four vertices $(j+\nu-1,j)$, $(j+\nu,j+1)$, $(j-1,j+\nu)$, and $(j, j+\nu+1)$. (Notice that the two entries $i,j$ of each vertex have been chosen in such a way that $h=i-j\ (\!\!\!\mod n)$ is not greater than $\nu$, as required.) These vertices  have respective vector entries $f_{\nu-1}\zeta^j$, $f_{\nu-1}\zeta^{j+1}$, $f_{\nu-1}\zeta^{j+\nu}$, and $f_{\nu-1}\zeta^{j+\nu+1}$. Then, the $\nu$-th entry of $\L\x=\lambda\x$ is 
 $$
 4f_{\nu}\zeta^j-f_{\nu-1}\zeta^j-f_{\nu-1}\zeta^{j+1}-f_{\nu-1}\zeta^{j+\nu}-f_{\nu-1}\zeta^{j+\nu-1}=\lambda f_{\nu}\zeta^j.
 $$ 
Thus,
$$
4f_{\nu}-f_{\nu-1}(1+\zeta+\zeta^{\nu}+\zeta^{\nu+1})=\lambda f_{\nu},
$$ 
which corresponds to the $\nu$-th entry of $\B(\zeta)\f=\lambda\f$.
 \item[$(iii)$] $h=\nu-1$: 
 The vertex $(j+\nu-1,j)$, with vector entry $x_{(j+\nu-1,j)}=f_{\nu-1}\zeta^j$, is adjacent to the four vertices 
 $(j+\nu-2,j)$, $(j+\nu,j)=(j,j+\nu)$, $(j+\nu-1,j-1)=(j-1,j+\nu-1)$, and $(j+\nu-1,j+1)$. To be consistent with the notation, these vertices  must have respective vector entries 
 \begin{equation}
 \label{eq:critic}
 f_{\nu-2}\zeta^j,\quad  f_{\nu}\zeta^{j}=f_{\nu}\zeta^{j+\nu},\quad f_{\nu}\zeta^{j-1}=f_{\nu}\zeta^{j-1+\nu},\quad \mbox{and}\quad f_{\nu-2}\zeta^{j+1}.
 \end{equation}
 Then, if the above equalities hold, the $(\nu-1)$-th entry of $\L\x=\lambda\x$ is 
 $$
 4f_{\nu-1}\zeta^j-f_{\nu-2}\zeta^j-f_{\nu}\zeta^{j}-f_{\nu}\zeta^{j-1}-f_{\nu-2}\zeta^{j+1}=\lambda f_{\nu-1}\zeta^j.
 $$ 
Thus,
$$
4f_{\nu-1}-f_{\nu-2}(1+\zeta)-f_{\nu}(1+\zeta^{-1})=\lambda  f_{\nu-1},
$$ 
which corresponds to the $(\nu-1)$-th entry of $\B(\zeta)\f=\lambda\f$.
Now, the second and third equalities in \eqref{eq:critic} hold
 if, either 
\begin{equation}
    f_{\nu}=0 \quad \mbox{ or } \quad \zeta^{\nu}=1.
    \label{eq:f=0oznu=1}
\end{equation}

Let us show that one or the other condition happens in the two following subcases:
\end{itemize}
\begin{itemize}
\item[$(iii.1)$] $r$ even:
When  $\zeta$ is an even power of $e^{i\frac{2\pi}{n}}$,
we have that $\zeta^{\nu}=1$, and \eqref{eq:f=0oznu=1} holds. Thus, each eigenvector of $\B^*(r)$ gives rise to an eigenvector of $\L$.
In particular, when $r=0$, the matrix  $\B^*(0)$ is similar to $\B(1)$, which 
equals the quotient matrix $\Q_L$ in \eqref{Qodd} of the path-shaped regular partition of $F_2(C_n)$. This is the special case when we know, in advance, that all eigenvalues of $\B^*(0)$ are eigenvalues of $\L$.
\item[$(iii.2)$] $r$ odd:
In this case $\cos(\frac{r\pi}{n})+\cos\left(\frac{r(n-1)\pi}{n}\right)=0$. Thus, the last row of $\B^*(r)$ is $(0,\ldots,0,4)$ and, hence, the matrix has one eigenvalue $\lambda=4$, with corresponding eigenvector $\f$ having $f_{\nu}\neq 0$. Moreover, since $\zeta$ is an odd power of $e^{i\frac{2\pi}{n}}$, we get $\zeta^{\nu}\neq 1$. Consequently, $\f$ does {\bf not} yield an eigenvector of $\L$ since  \eqref{eq:f=0oznu=1} does not hold.
Apart from this eigenvalue $4$, the other eigenvalues of $\B^*(r)$ are those of the principal submatrix $\B^-$ of the first $\nu-1$ rows and columns. Thus, the corresponding eigenvectors 
of $\B^*(r)$ have the last component $f_{\nu}=0$ to satisfy \eqref{eq:f=0oznu=1} and, then they yield an eigenvector of $\L$.
\end{itemize}
From the above discussion, we see that the eigenvalues not present in $\L$ are all equal to $4$ and belong to the spectrum of $\B^*(r)$ when $r$ is odd.
More precisely, for every odd $r(\neq n/2)$, there is an eigenvalue $4$ not in the spectrum of $\L$, obtaining $n/2$ of such eigenvalues. Besides, when $r=\nu=n/2$, all the off-diagonal entries of $\B^*(r)$ are zero, and $\B^*(r)=\diag(2,4,\ldots,4)$. 
If $\nu$ is even, this provides one eigenvalue $2$ and $n/2-1$ eigenvalues $4$ in $\L$. 
Otherwise, if $\nu$ is odd, the condition in \eqref{eq:critic} is not fulfilled for one $4$-eigenvector of $\B^*(r)$ and, hence, we only obtain one eigenvalue $2$ and $n/2-2$ eigenvalues $4$.
We obtain, in this
way, a total of $n\nu=2\nu^2$ eigenvalues (including repetitions), which is the number $2\nu(\nu-1)$ of eigenvalues
of the matrix $\L$ plus $\nu=|\Lambda|$ eigenvalues equal to $4$.
Now, to complete the proof, we need to show that the eigenvalues of $\B(z)$ for $z\in R(n)$ (not belonging to $\Lambda$) constitute the spectrum of $\L$. 
With this aim, we first use that, because of the properties of the polynomial matrix (see Dalf\'o, Fiol, Miller, Ryan, and \v{S}ir\'a\v{n} \cite{dfmrs17}), if 
$(\B(z)^{\ell})_{ii}=\alpha_{i0}^{(\ell)}+\alpha_{i1}^{(\ell)}z+\alpha_{i2}^{(\ell)}z^2+\cdots$ (for $\ell\geq0$), then 
$$
\tr(\L^{\ell})+\nu 4^{\ell}=\sum_{\lambda\in\spec \L\cup\Lambda}\lambda^{\ell} 
=n\sum_{i=1}^{\nu}\alpha_{i0}^{(\ell)},
$$
where we have taken into account that the matrices $\B(z)$, for $z\in R(n)$, have $\nu$ eigenvalues $4$ not in $\L$.
Since $\sum_{z\in R(n)}z^{\ell}=0$ for every $z\neq 1$
and $\ell\neq 0$, we have that
$$
\alpha_{i0}^{(\ell)}=\frac{1}{n}\sum_{z\in R(n)}(\B(z)^{\ell})_{ii}.
$$
Hence,
\begin{align*}
\sum_{\lambda\in\spec \L}\lambda^{\ell} +\nu 4^{\ell} &=\sum_{i=1}^{\nu}\sum_{z\in R(n)}(\B(z)^{\ell})_{ii}=
\sum_{z\in R(n)}\tr(\B(z)^{\ell})\\
 &=\sum_{z\in R(n)}
\sum_{\mu}\in\spec \B(z)\mu^\ell. 
\end{align*}
Since this holds for every $\ell\ge 0$, both multisets of eigenvalues in \eqref{spB(z)-spL} must coincide (see Gould \cite{go99}).
This completes the proof.
\end{proof}
By way of example, when $\L$ and $\B(z)$ are the matrices associated to $F_2(C_8)$, the equality $\tr(\L^{\ell})+\nu 4^{\ell}=\sum_{z\in R(n)}\tr(\B(z)^{\ell})$ for $\ell=0,\ldots,8$, give the values
32, 112, 512, 2656, 14976, 9792, 564032, 3670464, so that the corresponding traces  $\tr(\L^{\ell})=\displaystyle\sum_{\lambda\in\spec \L} \lambda^{\ell}$ are
28, 96, 448, 2400, 13952, 85696, 547648, 3604928,
as can be checked by using the values in \eqref{spF2(C8)} or Table \ref{taula:C8}.

\begin{table}[t!]
\begin{center}
\begin{tabular}{|c|cccc| }
\hline
$\zeta=e^{i\frac{2\pi}{8}}$, $z=\zeta^r$ & $\lambda_{r,1}$  & $\lambda_{r,2}$  & $\lambda_{r,3}$  & $\lambda_{r,4}$   \\
\hline\hline
$\spec(\B(\zeta^0))$ & \bf 0 & 1.506040792 & 4.890083735 &  7.603875471   \\
\hline
$\spec(\B(\zeta^1))=\spec(\B(\zeta^7))$ & \bf 0.5857864376 & 3.12596795 & 4.0 &  6.288245611   \\
\hline
$\spec(\B(\zeta^2))=\spec(\B(\zeta^6))$ & 0.9486257582 & \bf 2.0  &  4.517304045 & 6.534070196  \\
\hline
$\spec(\B(\zeta^3))=\spec(\B(\zeta^5))$ &  1.711754388  &  \bf 3.414213562 &   4.0 & 4.87403204  \\
\hline
$\spec(\B(\zeta^4))$ & 2.0   & \bf 4.0  &  4.0 &   4.0 \\
\hline
\end{tabular}
\end{center}
\caption{All the eigenvalues of the matrices $\B(\zeta^r)$, which yield the eigenvalues of the 2-token graph $F_2(C_8)$. The values in boldface correspond to the eigenvalues of $C_8$.}
\label{taula:C8}
\end{table}

\section{Asymptotic results}
In this last section, we derive closed formulas that give asymptotic approximations of the eigenvalues of $F_2(C_n)$
when $n$ is large.

\begin{theorem}
\label{th:cycles}
\begin{itemize}
 \item[$(i.1)$] 
 For $n$ odd, $n=2\nu+1$, and fixed odd $r<n$, the eigenvalues of $F_2(C_n)$, in the matrix $\B^*(r)$ of \eqref{B-star}, are asymptotically equal to
 \begin{equation}
 \label{asymp:odd-odd}
 \lambda_k=4+4\cos\left(\frac{r\pi}{n}\right)\cos\left(\frac{2k-1}{n-1}\pi\right),\quad  k=1,2,\ldots,\nu.
 \end{equation}
 \item[$(i.2)$] 
  For $n$ odd, $n=2\nu+1$, and fixed even $r<n$, the eigenvalues of the matrix $\B^*(r)$ in \eqref{B-star} are asymptotically equal to
 \begin{equation}
  \label{asymp:odd-even}
 \lambda_k=4+4\cos\left(\frac{r\pi}{n}\right)\cos\left(\frac{k-1}{n-1}2\pi\right)\quad k=1,2,\ldots,\nu.
 \end{equation}
 \item[$(ii)$] 
For $n$ even, $n=2\nu$, and fixed odd $r<n$ or $r=\nu$ even, the eigenvalues of $F_2(C_n)$, in the matrix $\B^*(r)$ of \eqref{B-star-even}, are asymptotically equal to
\begin{equation}
\label{asymp:even}
\lambda_k=4+4\cos\left(\frac{r\pi}{n}\right)\cos\left(\frac{2k-1}{n-1}\pi\right),\quad k=1,2,\ldots,\nu-1.
 \end{equation}
\end{itemize}
\end{theorem}

\begin{proof}
We begin with the case $(ii)$. 
For odd $r$ or $r=\nu$ even, the last row of  $\B^*(r)$ is $(0,\ldots,0,4)$. Hence, apart from the eigenvalue 4, the other eigenvalues of $\B^*(r)$ are those of the principal submatrix $\B^-$ of the first $\nu-1$ rows and columns. Moreover, since the function $\cos(\frac{r\pi}{x})$ is continuous for $x>1$ and tend to $1$ when $x\rightarrow \infty$,
we can work with the approximation 
\begin{align}
 \B^-&\approx \C_e=
\left(
\begin{array}{cccccc}
4-2\cos(\frac{r\pi}{n}) & 2\cos(\frac{r\pi}{n})& 0 & 0  & \ldots & 0\\
2\cos(\frac{r\pi}{n}) & 4 & 2\cos(\frac{r\pi}{n})& 0  & \ldots & 0\\
0 & 2\cos(\frac{r\pi}{n}) & 4 & 2\cos(\frac{r\pi}{n}) & \ldots & 0\\
\vdots  & \vdots & \ddots & \ddots &  \ddots & \vdots\\ 
0 & 0 &  \ldots & 2\cos(\frac{r\pi}{n}) &  4 & 2\cos(\frac{r\pi}{n})\\
0 & 0 & \ldots & 0 & 2\cos(\frac{r\pi}{n}) & 4 
\end{array}
\right).
\label{B-star-even-matrix}
\end{align}
Then, from the results by Yueh, see \cite[Th.2]{y05}, the eigenvalues of $\C_e$ are those in \eqref{asymp:even}.
 \\
For the cases $(i.1)$ and $(i.2)$ of odd $n=2\nu+1$, we proceed similarly. Now, the approximation  $\B'$ of the matrix $\B^*(r)$ in \eqref{B-star} is
\begin{align}
 \B'&\approx \C_o=
\left(
\begin{array}{cccccc}
4-2\cos(\frac{r\pi}{n}) & 2\cos(\frac{r\pi}{n})& 0 & 0  & \ldots & 0\\
2\cos(\frac{r\pi}{n}) & 4 & 2\cos(\frac{r\pi}{n})& 0  & \ldots & 0\\
0 & 2\cos(\frac{r\pi}{n}) & 4 & 2\cos(\frac{r\pi}{n}) & \ldots & 0\\
\vdots  & \vdots & \ddots & \ddots &  \ddots & \vdots\\ 
0 & 0 &  \ldots & 2\cos(\frac{r\pi}{n}) &  4 & 2\cos(\frac{r\pi}{n})\\
0 & 0 & \ldots & 0 & 2\cos(\frac{r\pi}{n}) & 4\pm 2\cos(\frac{r\pi}{n}) 
\end{array}
\right),
\label{B-star-odd-matrix}
\end{align}
where in $(\nu,\nu)$-entry, we must take the plus sign when $r$ is odd and the minus when $r$ is even.  
Then, by using the results of Yueh, \cite[Th.3]{y05} and
\cite[Th.5]{y05}, respectively, the eigenvalues of $\C_o$ are the claimed ones in \eqref{asymp:odd-odd} and \eqref{asymp:odd-even}.
\end{proof}
By using Gershgorin circles as in Proposition \ref{propo:B*(r)}, we can prove that the minimum eigenvalue of $F_2(C_n)$ coincides with the minimum eigenvalue of the matrices $\B^*(r)$ in \eqref{B-star} and \eqref{B-star-even} when $r=1$.
Moreover, the minimum eigenvalue in \eqref{asymp:odd-odd}
and \eqref{asymp:even} is  obtained when $k=\nu$ and $k=\nu-1$,
respectively, giving
\begin{equation}
\label{aprox-alpha-odd}
\alpha(F_2(C_n))\approx\lambda_{\nu}=4+4\cos\left(\frac{\pi}{n}\right)\cos\left(\frac{n-2}{n-1}\pi\right),
\end{equation}
and 
\begin{equation}
\label{aprox-alpha-even}
\alpha(F_2(C_n))\approx \lambda_{\nu-1}=4+4\cos\left(\frac{\pi}{n}\right)\cos\left(\frac{n-3}{n-1}\pi\right).
\end{equation}
Notice that, as $n$ increases, the expressions in \eqref{aprox-alpha-odd} and \eqref{aprox-alpha-even}  tend to $2-2\cos(\frac{2\pi}{n})$, which is the exact value of $\alpha(C_n)$, as expected.

Furthermore,  the equalities in \eqref{asymp:odd-even} for  $r=0$ become
\begin{equation}
  \label{asymp:r=0}
 \lambda_k=4+4\cos\left(\frac{k-1}{n-1}2\pi\right)=
 8\cos^2\left(\frac{k-1}{n-1}\pi\right),\quad k=1,2,\ldots,\nu,
 \end{equation}
which correspond to the asymptotic approximation of the exact values \eqref{lr-odd} in Proposition \ref{propo:path-shaped}$(ii)$.



\section*{Acknowledgements}

We thank Ruy Fabila-Monroy from Cinvestav (Mexico) and the anonymous referee because both of them helped us to make a better paper from its first version.


\end{document}